\documentclass[preprint,3p,a4paper,11pt,authoryear]{elsarticle}

\usepackage{setspace}

\makeatletter
\def\ps@pprintTitle{%
   \let\@oddhead\@empty
   \let\@oddfoot\@empty
   \let\@evenhead\@empty
   \let\@evenfoot\@empty
}
\makeatother

\usepackage[colorlinks=true,citecolor=blue]{hyperref}
\usepackage[round]{natbib}
\makeatletter

\patchcmd{\NAT@citex}
{\@citea\NAT@hyper@{%
		\NAT@nmfmt{\NAT@nm}%
		\hyper@natlinkbreak{\NAT@aysep\NAT@spacechar}{\@citeb\@extra@b@citeb}%
		\NAT@date}}
{\@citea\NAT@nmfmt{\NAT@nm}%
	\NAT@aysep\NAT@spacechar\NAT@hyper@{\NAT@date}}{}{}

\patchcmd{\NAT@citex}
{\@citea\NAT@hyper@{%
		\NAT@nmfmt{\NAT@nm}%
		\hyper@natlinkbreak{\NAT@spacechar\NAT@@open\if*#1*\else#1\NAT@spacechar\fi}%
		{\@citeb\@extra@b@citeb}%
		\NAT@date}}
{\@citea\NAT@nmfmt{\NAT@nm}%
	\NAT@spacechar\NAT@@open\if*#1*\else#1\NAT@spacechar\fi\NAT@hyper@{\NAT@date}}
{}{}

\makeatother

\journal{European Journal of Operational Research}

\biboptions{authoryear,sort&compress,round,comma}

\usepackage{amsthm}

\newtheorem{proposition}{Proposition}

\newtheorem{corollary}{Corollary}
\newtheorem{open}{Open Problem}

\usepackage{dsfont} 	
\usepackage{booktabs}
\usepackage{amsmath} 	
\usepackage{amssymb} 	

\makeatletter
\let\cl@chapter\undefined
\makeatletter
\usepackage{cleveref}
\crefname{appendix}{}{}

\crefname{condition}{Condition}{Conditions}
\crefname{section}{Section}{Sections}
\crefname{open}{Open Problem}{Open Problems}

\usepackage{pdfpages}

\usepackage{subcaption}
\usepackage{pgfplots}
\pgfplotsset{compat=1.17}

\usepackage{tikz}
\usetikzlibrary{shapes, arrows, positioning, matrix, calc, backgrounds, positioning, patterns, graphs, graphs.standard, math}

\definecolor{purple}{rgb}{0.6314,0.5412,0.8941}
\tikzstyle{block}=[rectangle, minimum width=2cm, text width=2cm, minimum height=0.5cm, text centered, draw, thick, rounded corners]
\tikzstyle{arrow}=[thick, ->, >=stealth, draw]

\tikzset{USA map/.cd,
state/.style={fill, draw=white, ultra thick},
HI/.style={}, AK/.style={}, FL/.style={}, NH/.style={}, MI/.style={}, MI/.style={}, VT/.style={}, ME/.style={}, RI/.style={}, NY/.style={}, PA/.style={}, NJ/.style={}, DE/.style={}, MD/.style={}, VA/.style={}, WV/.style={}, OH/.style={}, IN/.style={}, IL/.style={}, CT/.style={}, WI/.style={}, NC/.style={}, DC/.style={}, MA/.style={}, TN/.style={}, AR/.style={}, MO/.style={}, GA/.style={}, SC/.style={}, KY/.style={}, AL/.style={}, LA/.style={}, MS/.style={}, IA/.style={}, MN/.style={}, OK/.style={}, TX/.style={}, NM/.style={}, KS/.style={}, NE/.style={}, SD/.style={}, ND/.style={}, WY/.style={}, MT/.style={}, CO/.style={}, ID/.style={}, UT/.style={}, AZ/.style={}, NV/.style={}, OR/.style={}, WA/.style={}, CA/.style={}}

\tikzset{
every state/.style={USA map/state/.style={#1}},
HI/.style={USA map/HI/.style={#1}}, AK/.style={USA map/AK/.style={#1}}, FL/.style={USA map/FL/.style={#1}}, NH/.style={USA map/NH/.style={#1}}, MI/.style={USA map/MI/.style={#1}}, VT/.style={USA map/VT/.style={#1}}, ME/.style={USA map/ME/.style={#1}}, RI/.style={USA map/RI/.style={#1}}, NY/.style={USA map/NY/.style={#1}}, PA/.style={USA map/PA/.style={#1}}, NJ/.style={USA map/NJ/.style={#1}}, DE/.style={USA map/DE/.style={#1}}, MD/.style={USA map/MD/.style={#1}}, VA/.style={USA map/VA/.style={#1}}, WV/.style={USA map/WV/.style={#1}}, OH/.style={USA map/OH/.style={#1}}, IN/.style={USA map/IN/.style={#1}}, IL/.style={USA map/IL/.style={#1}}, CT/.style={USA map/CT/.style={#1}}, WI/.style={USA map/WI/.style={#1}}, NC/.style={USA map/NC/.style={#1}}, DC/.style={USA map/DC/.style={#1}}, MA/.style={USA map/MA/.style={#1}}, TN/.style={USA map/TN/.style={#1}}, AR/.style={USA map/AR/.style={#1}}, MO/.style={USA map/MO/.style={#1}}, GA/.style={USA map/GA/.style={#1}}, SC/.style={USA map/SC/.style={#1}}, KY/.style={USA map/KY/.style={#1}}, AL/.style={USA map/AL/.style={#1}}, LA/.style={USA map/LA/.style={#1}}, MS/.style={USA map/MS/.style={#1}}, IA/.style={USA map/IA/.style={#1}}, MN/.style={USA map/MN/.style={#1}}, OK/.style={USA map/OK/.style={#1}}, TX/.style={USA map/TX/.style={#1}}, NM/.style={USA map/NM/.style={#1}}, KS/.style={USA map/KS/.style={#1}}, NE/.style={USA map/NE/.style={#1}}, SD/.style={USA map/SD/.style={#1}}, ND/.style={USA map/ND/.style={#1}}, WY/.style={USA map/WY/.style={#1}}, MT/.style={USA map/MT/.style={#1}}, CO/.style={USA map/CO/.style={#1}}, ID/.style={USA map/ID/.style={#1}}, UT/.style={USA map/UT/.style={#1}}, AZ/.style={USA map/AZ/.style={#1}}, NV/.style={USA map/NV/.style={#1}}, OR/.style={USA map/OR/.style={#1}}, WA/.style={USA map/WA/.style={#1}}, CA/.style={USA map/CA/.style={#1}}
}

\newcommand{\USA}[1][]{
    \begin{scope}[y=0.80pt,x=0.80pt,yscale=-1, inner sep=0pt, outer sep=0pt,
    #1
    ]

    \path[USA map/state, USA map/FL, local bounding box=FL] (759.8167,439.1428) -- (762.0824,446.4614) -- (765.8121,456.2037) --
      (771.1468,465.5800) -- (774.8650,471.8847) -- (779.7149,477.3812) --
      (783.7564,481.0994) -- (785.3730,484.0093) -- (784.2414,485.3025) --
      (783.4330,486.5958) -- (786.3429,494.0322) -- (789.2528,496.9421) --
      (791.8394,502.2769) -- (795.3959,508.0967) -- (799.9224,516.3413) --
      (801.2157,523.9394) -- (801.7007,535.9023) -- (802.3473,537.6805) --
      (802.0240,541.0754) -- (799.5991,542.3687) -- (799.9224,544.3086) --
      (799.2758,546.2485) -- (799.5991,548.6734) -- (800.0841,550.6134) --
      (797.3358,553.8466) -- (794.2643,555.3015) -- (790.3844,555.4632) --
      (788.9295,557.0798) -- (786.5046,558.0497) -- (785.2113,557.5648) --
      (784.0797,556.5948) -- (783.7564,553.6849) -- (782.9481,550.2900) --
      (779.5532,545.1169) -- (775.9967,542.8537) -- (772.1168,542.5303) --
      (771.3085,543.8236) -- (768.2370,539.4588) -- (767.5903,535.9023) --
      (765.0037,531.8608) -- (763.2255,530.7291) -- (761.6089,532.8307) --
      (759.8306,532.5074) -- (757.7290,527.4959) -- (754.8191,523.6161) --
      (751.9092,518.2813) -- (749.3227,515.2097) -- (745.7662,511.4915) --
      (747.8677,509.0666) -- (751.1009,503.5702) -- (750.9393,501.9536) --
      (746.4128,500.9836) -- (744.7962,501.6302) -- (745.1195,502.2769) --
      (747.7061,503.2468) -- (746.2511,507.7733) -- (745.4428,508.2583) --
      (743.6646,504.2168) -- (742.3713,499.3670) -- (742.0480,496.6188) --
      (743.5029,491.9306) -- (743.5029,482.3927) -- (740.4314,478.6745) --
      (739.1381,475.6029) -- (733.9649,474.3096) -- (732.0250,473.6630) --
      (730.4084,471.0764) -- (727.0135,469.4598) -- (725.8819,466.0649) --
      (723.1337,465.0950) -- (720.7088,461.3768) -- (716.5056,459.9219) --
      (713.5957,458.4669) -- (711.0092,458.4669) -- (706.9676,459.2752) --
      (706.8060,461.2151) -- (707.6143,462.1851) -- (707.1293,463.3167) --
      (704.0578,463.1551) -- (700.3396,466.7116) -- (696.7830,468.6515) --
      (692.9032,468.6515) -- (689.6700,469.9448) -- (689.3466,467.1966) --
      (687.7300,465.2566) -- (684.8202,464.1250) -- (683.2036,462.6701) --
      (675.1205,458.7902) -- (667.5225,457.0120) -- (663.1577,457.6586) --
      (657.1762,458.1436) -- (651.1948,460.2452) -- (647.7155,460.8581) --
      (647.4776,452.8084) -- (644.8910,450.8685) -- (643.1128,449.0902) --
      (643.4361,446.0186) -- (653.6207,444.7254) -- (679.1631,441.8155) --
      (685.9529,441.1688) -- (691.3889,441.4491) -- (693.9754,445.3290) --
      (695.4304,446.7839) -- (703.5285,447.2991) -- (714.3483,446.6525) --
      (735.8607,445.3592) -- (741.3064,444.6848) -- (746.4140,444.8893) --
      (746.8408,447.7992) -- (749.0738,448.6075) -- (749.3087,443.9775) --
      (747.7805,439.8046) -- (749.0889,438.3647) -- (754.6436,438.8195) --
      (759.8167,439.1428) -- cycle(772.3621,571.5479) -- (774.7870,570.9012) --
      (776.0803,570.6588) -- (777.5353,568.3147) -- (779.8793,566.6980) --
      (781.1726,567.1830) -- (782.8701,567.5064) -- (783.2742,568.5571) --
      (779.7985,569.7696) -- (775.5953,571.2246) -- (773.2512,572.4370) --
      (772.3621,571.5479) -- cycle(785.8608,566.5364) -- (787.0733,567.5872) --
      (789.8215,565.4856) -- (795.1563,561.2824) -- (798.8745,557.4025) --
      (801.3803,550.7744) -- (802.3502,549.0770) -- (802.5119,545.6821) --
      (801.7844,546.1671) -- (800.8145,548.9962) -- (799.3595,553.6035) --
      (796.1263,558.8575) -- (791.7614,563.0607) -- (788.3666,565.0006) --
      (785.8608,566.5364) -- cycle;

    \path[USA map/state, USA map/NH, local bounding box=NH] (880.7990,142.4248) -- (881.6680,141.3483) -- (882.7582,138.0572) --
      (880.2152,137.1438) -- (879.7302,134.0722) -- (875.8503,132.9406) --
      (875.5270,130.1923) -- (868.2523,106.7515) -- (863.6508,92.2085) --
      (862.7538,92.2034) -- (862.1071,93.8200) -- (861.4605,93.3351) --
      (860.4905,92.3651) -- (859.0356,94.3050) -- (858.9871,99.3371) --
      (859.2987,105.0043) -- (861.2386,107.7525) -- (861.2386,111.7941) --
      (857.5204,116.8568) -- (854.9339,117.9885) -- (854.9339,119.1201) --
      (856.0655,120.8984) -- (856.0655,129.4664) -- (855.2572,138.6811) --
      (855.0955,143.5309) -- (856.0655,144.8242) -- (855.9038,149.3507) --
      (855.4188,151.1289) -- (856.3876,151.8382) -- (873.1753,147.4136) --
      (875.3502,146.8112) -- (877.1938,144.0378) -- (880.7990,142.4247) -- cycle;

    \begin{scope}
    \path[USA map/state, USA map/MI, local bounding box=MI] (697.8601,177.2369) -- (694.6269,168.9922)
        -- (692.3636,159.9392) -- (689.9387,156.7060) -- (687.3521,154.9277) --
        (685.7355,156.0594) -- (681.8557,157.8376) -- (679.9158,162.8491) --
        (677.1675,166.5673) -- (676.0359,167.2139) -- (674.5810,166.5673) .. controls
        (674.5810,166.5673) and (671.9944,165.1123) .. (672.1561,164.4657) .. controls
        (672.3177,163.8191) and (672.6410,159.4542) .. (672.6410,159.4542) --
        (676.0359,158.1609) -- (676.8442,154.7661) -- (677.4908,152.1795) --
        (679.9158,150.5629) -- (679.5924,140.5400) -- (677.9758,138.2767) --
        (676.6825,137.4684) -- (675.8742,135.3668) -- (676.6825,134.5585) --
        (678.2991,134.8818) -- (678.4608,133.2652) -- (676.0359,131.0020) --
        (674.7426,128.4154) -- (672.1561,128.4154) -- (667.6296,126.9605) --
        (662.1331,123.5656) -- (659.3849,123.5656) -- (658.7382,124.2123) --
        (657.7683,123.7273) -- (654.6967,121.4640) -- (651.7868,123.2423) --
        (648.8769,125.5055) -- (649.2003,129.0621) -- (650.1702,129.3854) --
        (652.2718,129.8704) -- (652.7568,130.6787) -- (650.1702,131.4870) --
        (647.5837,131.8103) -- (646.1287,133.5886) -- (645.8054,135.6901) --
        (646.1287,137.3067) -- (646.4520,142.8032) -- (642.8955,144.9048) --
        (642.2489,144.7431) -- (642.2489,140.5400) -- (643.5421,138.1151) --
        (644.1888,135.6901) -- (643.3805,134.8818) -- (641.4406,135.6901) --
        (640.4706,139.8933) -- (637.7224,141.0249) -- (635.9441,142.9649) --
        (635.7824,143.9348) -- (636.4291,144.7431) -- (635.7824,147.3297) --
        (633.5192,147.8147) -- (633.5192,148.9463) -- (634.3275,151.3712) --
        (633.1959,157.5143) -- (631.5793,161.5558) -- (632.2259,166.2440) --
        (632.7109,167.3756) -- (631.9026,169.8005) -- (631.5793,170.6088) --
        (631.2560,173.3570) -- (634.8125,179.3385) -- (637.7224,185.8049) --
        (639.1773,190.6547) -- (638.3690,195.3429) -- (637.3991,201.3243) --
        (634.9741,206.4974) -- (634.6508,209.2457) -- (631.3920,212.3308) --
        (635.8006,212.1688) -- (657.2191,209.9055) -- (664.4969,208.9184) --
        (664.5933,210.5848) -- (671.4452,209.3723) -- (681.7433,207.8692) --
        (685.5975,207.4083) -- (685.7356,206.8207) -- (685.8972,205.3658) --
        (687.9988,201.6476) -- (689.9994,199.9098) -- (689.7771,194.8579) --
        (691.3741,193.2609) -- (692.4647,192.9179) -- (692.6870,189.3614) --
        (694.2227,186.3303) -- (695.2735,186.9365) -- (695.4352,187.5832) --
        (696.2435,187.7448) -- (698.1834,186.7749) -- (697.8601,177.2369) -- cycle;

      \path[USA map/state, USA map/MI, local bounding box=MI2] (581.6193,82.0590) -- (583.4483,80.0014) --
        (585.6202,79.2012) -- (590.9929,75.3146) -- (593.2791,74.7431) --
        (593.7363,75.2003) -- (588.5923,80.3443) -- (585.2773,82.2876) --
        (583.2197,83.2021) -- (581.6193,82.0590) -- cycle(667.7937,114.1872) --
        (668.4403,116.6929) -- (671.6736,116.8546) -- (672.9668,115.6421) .. controls
        (672.9668,115.6421) and (672.8860,114.1872) .. (672.5627,114.0255) .. controls
        (672.2394,113.8639) and (670.9461,112.1664) .. (670.9461,112.1664) --
        (668.7637,112.4089) -- (667.1470,112.5706) -- (666.8237,113.7022) --
        (667.7937,114.1872) -- cycle(567.4921,111.2132) -- (568.2084,110.6328) --
        (570.9566,109.8245) -- (574.5131,107.5612) -- (574.5131,106.5913) --
        (575.1598,105.9446) -- (581.1412,104.9747) -- (583.5661,103.0347) --
        (587.9310,100.9331) -- (588.0926,99.6399) -- (590.0325,96.7300) --
        (591.8108,95.9217) -- (593.1041,94.1434) -- (595.3673,91.8802) --
        (599.7322,89.4553) -- (604.4203,88.9703) -- (605.5519,90.1019) --
        (605.2286,91.0719) -- (601.5104,92.0418) -- (600.0555,95.1134) --
        (597.7922,95.9217) -- (597.3073,98.3466) -- (594.8824,101.5798) --
        (594.5590,104.1664) -- (595.3673,104.6513) -- (596.3373,103.5197) --
        (599.8938,100.6098) -- (601.1871,101.9031) -- (603.4504,101.9031) --
        (606.6836,102.8731) -- (608.1385,104.0047) -- (609.5934,107.0762) --
        (612.3417,109.8245) -- (616.2215,109.6628) -- (617.6765,108.6928) --
        (619.2931,109.9861) -- (620.9097,110.4711) -- (622.2030,109.6628) --
        (623.3346,109.6628) -- (624.9512,108.6928) -- (628.9927,105.1363) --
        (632.3876,104.0047) -- (639.0157,103.6814) -- (643.5421,101.7414) --
        (646.1287,100.4482) -- (647.5837,100.6098) -- (647.5837,106.2679) --
        (648.0687,106.5913) -- (650.9785,107.3996) -- (652.9185,106.9146) --
        (659.0616,105.2980) -- (660.1932,104.1664) -- (661.6481,104.6513) --
        (661.6481,111.6027) -- (664.8813,114.6743) -- (666.1746,115.3209) --
        (667.4679,116.2909) -- (666.1746,116.6142) -- (665.3663,116.2909) --
        (661.6481,115.8059) -- (659.5465,116.4526) -- (657.2833,116.2909) --
        (654.0501,117.7458) -- (652.2718,117.7458) -- (646.4520,116.4526) --
        (641.2789,116.6142) -- (639.3390,119.2008) -- (632.3876,119.8474) --
        (629.9627,120.6557) -- (628.8311,123.7273) -- (627.5378,124.8589) --
        (627.0528,124.6972) -- (625.5978,123.0806) -- (621.0714,125.5055) --
        (620.4247,125.5055) -- (619.2931,123.8889) -- (618.4848,124.0506) --
        (616.5449,128.4154) -- (615.5749,132.4569) -- (612.3938,139.4577) --
        (611.2170,138.4235) -- (609.8453,137.3922) -- (607.9045,127.1041) --
        (604.3600,125.7341) -- (602.3074,123.4479) -- (590.1871,120.7044) --
        (587.3318,119.6747) -- (579.1014,117.5020) -- (571.2114,116.3589) --
        (567.4921,111.2132) -- cycle;

    \end{scope}
    \path[USA map/state, USA map/VT, local bounding box=VT] (844.4842,154.0579) -- (844.8009,148.7123) -- (841.9101,137.9281) --
      (841.2635,137.6048) -- (838.3536,136.3115) -- (839.1619,133.4016) --
      (838.3536,131.3000) -- (835.6536,126.6600) -- (836.6235,122.7802) --
      (835.8152,117.6070) -- (833.3903,111.1406) -- (832.5847,106.2181) --
      (859.0041,99.4863) -- (859.3128,105.0085) -- (861.2291,107.7507) --
      (861.2291,111.7922) -- (857.5219,116.8502) -- (854.9353,117.9929) --
      (854.9243,119.1135) -- (856.2343,120.6326) -- (855.9234,128.7305) --
      (855.3139,137.9894) -- (855.0860,143.5463) -- (856.0560,144.8396) --
      (855.8943,149.4103) -- (855.4093,151.1002) -- (856.4235,151.8274) --
      (848.9860,153.3341) -- (844.4842,154.0579) -- cycle;

    \path[USA map/state, USA map/ME, local bounding box=ME] (922.8398,78.8307) -- (924.7797,80.9323) -- (927.0429,84.6505) --
      (927.0429,86.5904) -- (924.9413,91.2786) -- (923.0014,91.9252) --
      (919.6065,94.9968) -- (914.7567,100.4932) .. controls (914.7567,100.4932) and
      (914.1101,100.4932) .. (913.4635,100.4932) .. controls (912.8168,100.4932) and
      (912.4935,98.3916) .. (912.4935,98.3916) -- (910.7152,98.5533) --
      (909.7453,100.0082) -- (907.3204,101.4632) -- (906.3504,102.9181) --
      (907.9670,104.3731) -- (907.4820,105.0197) -- (906.9970,107.7679) --
      (905.0571,107.6063) -- (905.0571,105.9897) -- (904.7338,104.6964) --
      (903.2789,105.0197) -- (901.5006,101.7865) -- (899.3990,103.0798) --
      (900.6923,104.5347) -- (901.0156,105.6664) -- (900.2073,106.9596) --
      (900.5306,110.0312) -- (900.6923,111.6478) -- (899.0757,114.2344) --
      (896.1658,114.7193) -- (895.8425,117.6292) -- (890.5077,120.7008) --
      (889.2144,121.1858) -- (887.5978,119.7308) -- (884.5262,123.2873) --
      (885.4962,126.5206) -- (884.0412,127.8138) -- (883.8796,132.1787) --
      (882.7563,138.4380) -- (880.2941,137.2821) -- (879.8091,134.2105) --
      (875.9292,133.0789) -- (875.6059,130.3306) -- (868.3311,106.8898) --
      (863.6326,92.2501) -- (865.0531,92.1319) -- (866.5669,92.5418) --
      (866.5669,89.9553) -- (867.8752,85.4588) -- (870.4618,80.7706) --
      (871.9167,76.7291) -- (869.9768,74.3042) -- (869.9768,68.3228) --
      (870.7851,67.3528) -- (871.5934,64.6046) -- (871.4317,63.1497) --
      (871.2701,58.2998) -- (873.0483,53.4500) -- (875.9582,44.5587) --
      (878.0598,40.3555) -- (879.3531,40.3555) -- (880.6464,40.5172) --
      (880.6464,41.6488) -- (881.9397,43.9121) -- (884.6879,44.5587) --
      (885.4962,43.7504) -- (885.4962,42.7804) -- (889.5377,39.8705) --
      (891.3160,38.0923) -- (892.7709,38.2539) -- (898.7523,40.6788) --
      (900.6923,41.6488) -- (909.7453,71.5560) -- (915.7267,71.5560) --
      (916.5350,73.4959) -- (916.6967,78.3457) -- (919.6066,80.6090) --
      (920.4149,80.6090) -- (920.5765,80.1240) -- (920.0915,78.9924) --
      (922.8398,78.8307) -- cycle(901.9080,108.9783) -- (903.4438,107.4425) --
      (904.8179,108.4933) -- (905.3837,110.9182) -- (903.6863,111.8073) --
      (901.9080,108.9782) -- cycle(908.6169,103.0776) -- (910.3952,104.9367) ..
      controls (910.3952,104.9367) and (911.6885,105.0175) .. (911.6885,104.6942) ..
      controls (911.6885,104.3709) and (911.9310,102.6735) .. (911.9310,102.6735) --
      (912.8201,101.8652) -- (912.0118,100.0869) -- (909.9911,100.8144) --
      (908.6169,103.0776) -- cycle;

    \path[USA map/state, USA map/RI, local bounding box=RI] (874.0700,178.8954) -- (870.3742,163.9394) -- (876.6435,162.0942) --
      (878.8346,164.0214) -- (882.1411,168.3420) -- (884.8290,172.7441) --
      (881.8297,174.3689) -- (880.5364,174.2072) -- (879.4048,175.9855) --
      (876.9799,177.9254) -- (874.0700,178.8954) -- cycle;

    \path[USA map/state, USA map/NY, local bounding box=NY] (830.3794,188.7456) -- (829.2478,187.7756) -- (826.6612,187.6140) --
      (824.3980,185.6741) -- (822.7674,179.5449) -- (819.3089,179.6354) --
      (816.8652,176.9272) -- (797.4799,181.3092) -- (754.4781,190.0389) --
      (746.9485,191.2669) -- (746.2103,184.7985) -- (747.6384,183.6731) --
      (748.9317,182.5415) -- (749.9017,180.9249) -- (751.6799,179.7933) --
      (753.6198,178.0150) -- (754.1048,176.3984) -- (756.2064,173.6502) --
      (757.3380,172.6802) -- (757.1764,171.7103) -- (755.8831,168.6387) --
      (754.1048,168.4770) -- (752.1649,162.3339) -- (755.0748,160.5557) --
      (759.4396,159.1007) -- (763.4811,157.8074) -- (766.7143,157.3225) --
      (773.0191,157.1608) -- (774.9590,158.4541) -- (776.5756,158.6158) --
      (778.6772,157.3225) -- (781.2638,156.1908) -- (786.4369,155.7059) --
      (788.5385,153.9276) -- (790.3168,150.6944) -- (791.9334,148.7545) --
      (794.0350,148.7545) -- (795.9749,147.6228) -- (796.1365,145.3596) --
      (794.6816,143.2580) -- (794.3583,141.8031) -- (795.4899,139.7015) --
      (795.4899,138.2465) -- (793.7116,138.2465) -- (791.9334,137.4382) --
      (791.1251,136.3066) -- (790.9634,133.7200) -- (796.7832,128.2236) --
      (797.4298,127.4153) -- (798.8848,124.5054) -- (801.7947,119.9789) --
      (804.5429,116.2607) -- (806.6445,113.8358) -- (809.0596,112.0102) --
      (812.1409,110.7643) -- (817.6374,109.4710) -- (820.8706,109.6326) --
      (825.3971,108.1777) -- (832.9623,106.1065) -- (833.4821,111.0862) --
      (835.9070,117.5526) -- (836.7153,122.7258) -- (835.7453,126.6056) --
      (838.3319,131.1321) -- (839.1402,133.2337) -- (838.3319,136.1436) --
      (841.2418,137.4369) -- (841.8884,137.7602) -- (844.9600,148.7532) --
      (844.4237,153.8128) -- (843.9387,164.6441) -- (844.7470,170.1406) --
      (845.5553,173.6971) -- (847.0103,180.9719) -- (847.0103,189.0549) --
      (845.8787,191.3182) -- (847.7180,193.3109) -- (848.5145,194.9894) --
      (846.5746,196.7676) -- (846.8979,198.0609) -- (848.1912,197.7376) --
      (849.6462,196.4443) -- (851.9094,193.8577) -- (853.0410,193.2111) --
      (854.6576,193.8577) -- (856.9209,194.0194) -- (864.8422,190.1396) --
      (867.7521,187.3913) -- (869.0454,185.9364) -- (873.2486,187.5530) --
      (869.8537,191.1095) -- (865.9739,194.0194) -- (858.8608,199.3542) --
      (856.2742,200.3242) -- (850.4545,202.2641) -- (846.4130,203.3957) --
      (845.2382,202.8628) -- (844.9942,199.1743) -- (845.4792,196.4260) --
      (845.3175,194.3244) -- (842.5040,192.6254) -- (837.9775,191.6555) --
      (834.0976,190.5238) -- (830.3794,188.7456) -- cycle;

    \path[USA map/state, USA map/PA, local bounding box=PA] (825.1237,224.6920) -- (826.4321,224.4211) -- (828.7616,223.1678) --
      (829.9735,220.6847) -- (831.5901,218.4215) -- (834.8233,215.3499) --
      (834.8233,214.5416) -- (832.3984,212.9250) -- (828.8419,210.5001) --
      (827.8719,207.9135) -- (825.1237,207.5902) -- (824.9620,206.4586) --
      (824.1537,203.7103) -- (826.4170,202.5787) -- (826.5787,200.1538) --
      (825.2854,198.8605) -- (825.4470,197.2439) -- (827.3870,194.1724) --
      (827.3870,191.1008) -- (830.0846,188.4549) -- (829.1643,187.7799) --
      (826.6402,187.5870) -- (824.3457,185.6471) -- (822.7958,179.5310) --
      (819.2912,179.6316) -- (816.8360,176.9282) -- (798.7450,181.1260) --
      (755.7432,189.8557) -- (746.8519,191.3106) -- (746.2312,184.7892) --
      (740.8687,189.8569) -- (739.5754,190.3419) -- (735.3731,193.3508) --
      (738.2839,212.4882) -- (740.7655,222.2176) -- (744.3373,241.4791) --
      (747.6066,240.8414) -- (759.5502,239.3389) -- (797.4768,231.6737) --
      (812.3531,228.8504) -- (820.6534,227.2280) -- (820.9205,226.9895) --
      (823.0221,225.3729) -- (825.1237,224.6920) -- cycle;

    \path[USA map/state, USA map/NJ, local bounding box=NJ] (829.6794,188.4602) -- (827.3569,191.1944) -- (827.3569,194.2660) --
      (825.4169,197.3375) -- (825.2553,198.9542) -- (826.5486,200.2474) --
      (826.3869,202.6724) -- (824.1237,203.8040) -- (824.9320,206.5522) --
      (825.0936,207.6838) -- (827.8419,208.0072) -- (828.8118,210.5937) --
      (832.3684,213.0187) -- (834.7933,214.6353) -- (834.7933,215.4436) --
      (831.8101,218.1401) -- (830.1934,220.4034) -- (828.7385,223.1516) --
      (826.4752,224.4449) -- (826.0128,226.0474) -- (825.7703,227.2598) --
      (825.1611,229.8666) -- (826.2533,232.1108) -- (829.4865,235.0206) --
      (834.3364,237.2839) -- (838.3779,237.9305) -- (838.5395,239.3855) --
      (837.7312,240.3554) -- (838.0545,243.1037) -- (838.8628,243.1037) --
      (840.9644,240.6788) -- (841.7727,235.8289) -- (844.5210,231.7874) --
      (847.5925,225.3210) -- (848.7241,219.8246) -- (848.0775,218.6929) --
      (847.9158,209.3166) -- (846.2992,205.9218) -- (845.1676,206.7301) --
      (842.4194,207.0534) -- (841.9344,206.5684) -- (843.0660,205.5984) --
      (845.1676,203.6585) -- (845.2307,202.5647) -- (844.8463,199.1308) --
      (845.4197,196.3826) -- (845.3022,194.4136) -- (842.4947,192.6632) --
      (837.4025,191.4875) -- (833.2651,190.1059) -- (829.6795,188.4602) -- cycle;

    \path[USA map/state, USA map/DE, local bounding box=DE] (825.6261,228.2791) -- (825.9944,226.1322) -- (826.3695,224.4412) --
      (824.7465,224.8389) -- (823.1310,225.3065) -- (820.9248,227.0708) --
      (822.6449,232.1137) -- (824.9081,237.7718) -- (827.0097,247.4714) --
      (828.6263,253.7762) -- (833.6378,253.6145) -- (839.7799,252.4339) --
      (837.5157,245.0476) -- (836.5457,245.5326) -- (832.9892,243.1077) --
      (831.2109,238.4195) -- (829.2710,234.8630) -- (826.1239,231.9927) --
      (825.2597,229.8946) -- (825.6261,228.2791) -- cycle;

    \path[USA map/state, USA map/MD, local bounding box=MD] (839.7917,252.4148) -- (833.7832,253.6186) -- (828.6403,253.7361) --
      (826.7967,246.8137) -- (824.8719,237.6444) -- (822.2993,231.4560) --
      (821.0109,227.0576) -- (813.5049,228.6800) -- (798.6287,231.5033) --
      (761.1773,239.0542) -- (762.3086,244.0659) -- (763.2785,249.7240) --
      (763.6018,249.4007) -- (765.7034,246.9758) -- (767.9667,244.3581) --
      (770.3916,243.7425) -- (771.8466,242.2876) -- (773.6248,239.7010) --
      (774.9181,240.3477) -- (777.8280,240.0243) -- (780.4146,237.9228) --
      (782.4215,236.4695) -- (784.2667,235.9845) -- (785.9110,237.1145) --
      (788.8209,238.5694) -- (790.7609,240.3477) -- (791.9733,241.8835) --
      (796.0957,243.5809) -- (796.0957,246.4908) -- (801.5921,247.7841) --
      (802.7366,248.3260) -- (804.1485,246.2977) -- (807.0304,248.2679) --
      (805.7523,250.7498) -- (804.9870,254.7355) -- (803.2087,257.3220) --
      (803.2087,259.4236) -- (803.8554,261.2019) -- (808.9193,262.5576) --
      (813.2304,262.4959) -- (816.3020,263.4659) -- (818.4035,263.7892) --
      (819.3735,261.6876) -- (817.9186,259.5860) -- (817.9186,257.8077) --
      (815.4937,255.7062) -- (813.3921,250.2097) -- (814.6854,244.8749) --
      (814.5237,242.7733) -- (813.2304,241.4800) .. controls (813.2304,241.4800) and
      (814.6854,239.8634) .. (814.6854,239.2168) .. controls (814.6854,238.5701) and
      (815.1703,237.1152) .. (815.1703,237.1152) -- (817.1103,235.8219) --
      (819.0502,234.2053) -- (819.5352,235.1753) -- (818.0802,236.7919) --
      (816.7869,240.5101) -- (817.1103,241.6417) -- (818.8885,241.9650) --
      (819.3735,247.4615) -- (817.2719,248.4314) -- (817.5952,251.9880) --
      (818.0802,251.8263) -- (819.2118,249.8864) -- (820.8285,251.6646) --
      (819.2118,252.9579) -- (818.8885,256.3528) -- (821.4751,259.7477) --
      (825.3549,260.2327) -- (826.9716,259.4244) -- (830.2081,263.6073) --
      (831.5665,264.1436) -- (838.2201,261.3466) -- (840.2277,257.3228) --
      (839.7917,252.4148) -- cycle(823.8222,261.4435) -- (824.9538,263.9492) --
      (825.1155,265.7275) -- (826.2471,267.5866) .. controls (826.2471,267.5866) and
      (827.1362,266.6975) .. (827.1362,266.3741) .. controls (827.1362,266.0508) and
      (826.4087,263.3026) .. (826.4087,263.3026) -- (825.6813,260.9585) --
      (823.8222,261.4435) -- cycle;

    \path[USA map/state, USA map/VA, local bounding box=VA] (831.6389,266.0689) -- (831.4949,264.1219) -- (837.9484,261.5720) --
      (837.1780,264.7899) -- (834.2580,268.5690) -- (833.8399,273.1548) --
      (834.3017,276.5452) -- (832.4737,281.5234) -- (830.3094,283.4395) --
      (828.8391,278.7987) -- (829.2850,273.3496) -- (830.8720,269.1665) --
      (831.6389,266.0689) -- cycle(834.9790,294.3703) -- (776.8049,306.9457) --
      (739.3779,312.2248) -- (732.6996,311.8496) -- (730.1143,313.7760) --
      (722.7752,313.9967) -- (714.3931,314.9743) -- (703.4781,316.5890) --
      (713.9475,310.9778) -- (713.9344,308.9028) -- (715.4545,306.7567) --
      (726.0083,295.2553) -- (729.9550,299.7327) -- (733.7380,300.6967) --
      (736.2815,299.5564) -- (738.5187,298.2452) -- (741.0553,299.5887) --
      (744.9695,298.1610) -- (746.8462,293.6046) -- (749.4471,294.1447) --
      (752.3024,292.0134) -- (754.1016,292.5070) -- (756.9288,288.8304) --
      (757.2771,286.7473) -- (756.3134,285.4718) -- (757.3162,283.6051) --
      (762.5905,271.3280) -- (763.2072,265.5929) -- (764.4361,265.0694) --
      (766.6147,267.5122) -- (770.5505,267.2111) -- (772.4797,259.6374) --
      (775.2737,259.0766) -- (776.3235,256.3355) -- (778.9033,253.9886) --
      (781.6751,248.2934) -- (781.7600,243.2259) -- (791.5815,247.0487) .. controls
      (792.2624,247.3891) and (792.4144,241.9996) .. (792.4144,241.9996) --
      (796.0670,243.5979) -- (796.1353,246.5361) -- (801.9195,247.8355) --
      (804.0525,249.0117) -- (805.7124,251.0674) -- (805.0578,254.7161) --
      (803.1104,257.3071) -- (803.2202,259.3662) -- (803.8092,261.2191) --
      (808.7880,262.4875) -- (813.2392,262.5274) -- (816.3081,263.4860) --
      (818.2516,263.7953) -- (818.9664,266.8838) -- (822.1568,267.2863) --
      (823.0249,268.4863) -- (822.5854,273.1764) -- (823.9601,274.2790) --
      (823.4812,276.2094) -- (824.7106,276.9991) -- (824.4888,278.3837) --
      (821.7948,278.2888) -- (821.8838,279.9044) -- (824.1648,281.4472) --
      (824.2863,282.8591) -- (826.0594,284.6445) -- (826.5512,287.1686) --
      (823.9982,288.5499) -- (825.5704,290.0442) -- (831.3714,288.3584) --
      (834.9790,294.3703) -- cycle;

    \path[USA map/state, USA map/WV, local bounding box=WV] (761.1855,238.9673) -- (762.2975,243.9118) -- (763.3810,249.9432) --
      (765.5113,247.3628) -- (767.7745,244.2913) -- (770.3129,243.6757) --
      (771.7678,242.2208) -- (773.5461,239.6342) -- (774.9911,240.2808) --
      (777.9010,239.9575) -- (780.4875,237.8559) -- (782.4944,236.4027) --
      (784.3397,235.9177) -- (785.6436,236.9342) -- (789.2868,238.7558) --
      (791.2268,240.5341) -- (792.6009,241.8273) -- (791.8392,247.3823) --
      (786.0042,244.8411) -- (781.7590,243.2190) -- (781.6579,248.3975) --
      (778.9102,253.9342) -- (776.3802,256.3609) -- (775.1881,259.1102) --
      (772.5445,259.6103) -- (771.6467,263.2122) -- (770.6034,267.1619) --
      (766.6352,267.5026) -- (764.3115,265.0638) -- (763.2403,265.6232) --
      (762.6076,271.0929) -- (761.2574,274.6274) -- (756.2990,285.5823) --
      (757.1956,286.7430) -- (756.9898,288.6516) -- (754.1811,292.5360) --
      (752.3726,291.9918) -- (749.4045,294.1515) -- (746.8622,293.5793) --
      (744.8629,298.1349) .. controls (744.8629,298.1349) and (741.6036,299.5651) ..
      (740.9400,299.5026) .. controls (740.7795,299.4875) and (738.4709,298.2535) ..
      (738.4709,298.2535) -- (736.1344,299.6329) -- (733.7246,300.6773) --
      (729.9799,299.7881) -- (728.8585,298.6199) -- (726.6663,295.5965) --
      (723.5237,293.6084) -- (721.8121,289.9851) -- (717.5273,286.5169) --
      (716.8806,284.2537) -- (714.2940,282.7987) -- (713.4857,281.1821) --
      (713.2432,275.9282) -- (715.4257,275.8474) -- (717.3656,275.0391) --
      (717.5273,272.2908) -- (719.1439,270.8359) -- (719.3055,265.8244) --
      (720.2755,261.9445) -- (721.5688,261.2979) -- (722.8620,262.4295) --
      (723.3470,264.2078) -- (725.1253,263.2378) -- (725.6103,261.6212) --
      (724.4787,259.8430) -- (724.4787,257.4180) -- (725.4486,256.1248) --
      (727.7119,252.7299) -- (729.0052,251.2749) -- (731.1068,251.7599) --
      (733.3700,250.1433) -- (736.4415,246.7484) -- (738.7048,242.8686) --
      (739.0281,237.2105) -- (739.5131,232.1990) -- (739.5131,227.5108) --
      (738.3815,224.4393) -- (739.3514,222.9843) -- (740.6349,221.6910) --
      (744.1262,241.5181) -- (748.7572,240.7670) -- (761.1855,238.9673) -- cycle;

    \path[USA map/state, USA map/OH, local bounding box=OH] (735.3250,193.3283) -- (729.2314,197.3817) -- (725.3516,199.6449) --
      (721.9567,203.3631) -- (717.9152,207.2430) -- (714.6820,208.0513) --
      (711.7721,208.5362) -- (706.2756,211.1228) -- (704.1741,211.2845) --
      (700.7792,208.2129) -- (695.6061,208.8596) -- (693.0195,207.4046) --
      (690.6384,206.0538) -- (685.7459,206.7572) -- (675.5612,208.3738) --
      (664.3544,210.5585) -- (665.6477,225.1888) -- (667.4259,238.9300) --
      (670.0125,262.3708) -- (670.5783,267.2020) -- (674.7007,267.0729) --
      (677.1256,266.2646) -- (680.4894,267.7678) -- (682.5598,272.1326) --
      (687.6988,272.1155) -- (689.5905,274.2342) -- (691.3517,274.1689) --
      (693.8901,272.8274) -- (696.3943,273.1989) -- (701.8155,273.6816) --
      (703.5425,271.5489) -- (705.8882,270.2557) -- (707.9587,269.5748) --
      (708.6053,272.3230) -- (710.3836,273.2930) -- (713.8593,275.6371) --
      (716.0417,275.5563) -- (717.3748,275.0638) -- (717.5595,272.3023) --
      (719.1449,270.8473) -- (719.2441,266.0546) .. controls (719.2441,266.0546) and
      (720.2680,261.9455) .. (720.2680,261.9455) -- (721.5673,261.3442) --
      (722.8887,262.4920) -- (723.4268,264.1890) -- (725.1459,263.1516) --
      (725.5849,261.6908) -- (724.4682,259.7878) -- (724.5345,257.4733) --
      (725.2835,256.4010) -- (727.4363,253.0946) -- (728.4865,251.5512) --
      (730.5881,252.0362) -- (732.8513,250.4196) -- (735.9229,247.0247) --
      (738.6944,242.9460) -- (739.0147,237.8905) -- (739.4997,232.8790) --
      (739.3229,227.5721) -- (738.3681,224.6773) -- (738.7193,223.4875) --
      (740.5237,221.7374) -- (738.2349,212.6901) -- (735.3250,193.3283) -- cycle;

    \path[USA map/state, USA map/IN, local bounding box=IN] (619.5695,299.9713) -- (619.6348,297.1127) -- (620.1198,292.5862) --
      (622.3831,289.6764) -- (624.1613,285.7965) -- (626.7479,281.5933) --
      (626.2629,275.7735) -- (624.4847,273.0253) -- (624.1613,269.7921) --
      (624.9697,264.2956) -- (624.4847,257.3442) -- (623.1914,241.3398) --
      (621.8981,225.9820) -- (620.9276,214.2620) -- (623.9987,215.1515) --
      (625.4536,216.1215) -- (626.5853,215.7982) -- (628.6868,213.8582) --
      (631.5164,212.2413) -- (636.6092,212.0792) -- (658.5951,209.8160) --
      (664.1708,209.2828) -- (665.6740,225.2390) -- (669.9253,262.0806) --
      (670.5238,267.8521) -- (670.1523,270.1154) -- (671.3802,271.9108) --
      (671.4766,273.2833) -- (668.9554,274.8828) -- (665.4159,276.4341) --
      (662.2138,276.9844) -- (661.6153,281.8514) -- (657.0406,285.1638) --
      (654.2442,289.1743) -- (654.5675,291.5510) -- (653.9862,293.0852) --
      (650.6597,293.0852) -- (649.0742,291.4686) -- (646.5809,292.7308) --
      (643.8979,294.2339) -- (644.0596,297.2884) -- (642.8658,297.5464) --
      (642.3979,296.5283) -- (640.2311,295.0251) -- (636.9807,296.3666) --
      (635.4294,299.3729) -- (633.9916,298.5646) -- (632.5366,296.9651) --
      (628.0723,297.4500) -- (622.4795,298.4200) -- (619.5696,299.9713) -- cycle;

    \path[USA map/state, USA map/IL, local bounding box=IL] (619.5415,300.3424) -- (619.5727,297.1127) -- (620.1400,292.4668) --
      (622.4726,289.5509) -- (624.3392,285.4751) -- (626.5722,281.4798) --
      (626.2007,276.2274) -- (624.1955,272.6848) -- (624.0991,269.3382) --
      (624.7940,264.0687) -- (623.9686,256.8903) -- (622.9022,241.1128) --
      (621.6089,226.0955) -- (620.6867,214.4563) -- (620.4141,213.5349) --
      (619.6058,210.9483) -- (618.3126,207.2301) -- (616.6960,205.4519) --
      (615.2410,202.8653) -- (615.0074,197.3764) -- (569.2110,199.9746) --
      (569.4396,202.3466) -- (571.7259,203.0324) -- (572.6403,204.1755) --
      (573.0976,206.0045) -- (576.9842,209.4339) -- (577.6701,211.7201) --
      (576.9842,215.1494) -- (575.1552,218.8074) -- (574.4693,221.3222) --
      (572.1831,223.1512) -- (570.3541,223.8371) -- (565.0958,225.2088) --
      (564.4099,227.0378) -- (563.7241,229.0954) -- (564.4099,230.4672) --
      (566.2389,232.0675) -- (566.0103,236.1827) -- (564.1813,237.7831) --
      (563.4954,239.3834) -- (563.4954,242.1269) -- (561.6665,242.5841) --
      (560.0661,243.7273) -- (559.8375,245.0990) -- (560.0661,247.1566) --
      (558.3514,248.4712) -- (557.3226,251.2718) -- (557.7799,254.9298) --
      (560.0661,262.2457) -- (567.3820,269.7902) -- (572.8690,273.4482) --
      (572.6403,277.7920) -- (573.5548,279.1638) -- (579.9563,279.6210) --
      (582.6997,280.9928) -- (582.0139,284.6507) -- (579.7277,290.5949) --
      (579.0418,293.7956) -- (581.3280,297.6822) -- (587.7294,302.9405) --
      (592.3019,303.6264) -- (594.3595,308.6561) -- (596.4171,311.8568) --
      (595.5026,314.8289) -- (597.1030,318.9441) -- (598.9319,321.0017) --
      (600.3460,320.1210) -- (601.2536,318.0462) -- (603.4667,316.2990) --
      (605.5982,315.6846) -- (608.2007,316.8644) -- (611.8277,318.2401) --
      (613.0167,317.9419) -- (613.2165,315.6834) -- (611.9292,313.2717) --
      (612.2334,310.8949) -- (614.0718,309.5475) -- (617.0944,308.7372) --
      (618.3553,308.2787) -- (617.7427,306.8918) -- (616.9513,304.5374) --
      (618.3839,303.5565) -- (619.5414,300.3424) -- cycle;

    \path[USA map/state, USA map/CT, local bounding box=CT] (874.0683,178.8629) -- (870.3909,163.9841) -- (865.6721,164.9044) --
      (844.4433,169.6475) -- (845.4435,172.8731) -- (846.8984,180.1479) --
      (847.0752,189.1148) -- (845.8552,191.2897) -- (847.7760,193.2220) --
      (852.0475,189.3164) -- (855.6040,186.0832) -- (857.5439,183.9816) --
      (858.3523,184.6282) -- (861.1005,183.1733) -- (866.2736,182.0417) --
      (874.0683,178.8629) -- cycle;

    \path[USA map/state, USA map/WI, local bounding box=WI] (615.0659,197.3687) -- (614.9992,194.2112) -- (613.8201,189.6847) --
      (613.1734,183.5417) -- (612.0418,181.1167) -- (613.0118,178.0452) --
      (613.8201,175.1353) -- (615.2750,172.5487) -- (614.6284,169.1539) --
      (613.9817,165.5973) -- (614.4667,163.8191) -- (616.4066,161.3942) --
      (616.5683,158.6459) -- (615.7600,157.3526) -- (616.4066,154.7661) --
      (615.9541,150.5954) -- (618.7024,144.9373) -- (621.6122,138.1475) --
      (621.7739,135.8843) -- (621.4506,134.9143) -- (620.6423,135.3993) --
      (616.4391,141.7041) -- (613.6909,145.7456) -- (611.7510,147.5238) --
      (610.9427,149.7871) -- (608.9877,150.5954) -- (607.8561,152.5353) --
      (606.4011,152.2120) -- (606.2395,150.4337) -- (607.5328,148.0088) --
      (609.6343,143.3207) -- (611.4126,141.7040) -- (612.4034,139.3462) --
      (609.8430,137.4449) -- (607.8682,127.0779) -- (604.3207,125.7359) --
      (602.3744,123.4276) -- (590.2447,120.7059) -- (587.3688,119.6939) --
      (579.1557,117.5266) -- (571.2378,116.3678) -- (567.4726,111.2372) --
      (566.7222,111.7912) -- (565.5243,111.6295) -- (564.8777,110.4979) --
      (563.5437,110.7944) -- (562.4120,110.9561) -- (560.6338,111.9261) --
      (559.6638,111.2794) -- (560.3105,109.3395) -- (562.2504,106.2679) --
      (563.3820,105.1363) -- (561.4421,103.6814) -- (559.3405,104.4897) --
      (556.4306,106.4296) -- (548.9942,109.6628) -- (546.0843,110.3094) --
      (543.1745,109.8245) -- (542.1927,108.9462) -- (540.0760,111.7814) --
      (539.8474,114.5249) -- (539.8474,122.9839) -- (538.7043,124.5843) --
      (533.4460,128.4708) -- (531.1597,134.4150) -- (531.6170,134.6437) --
      (534.1318,136.7013) -- (534.8177,139.9020) -- (532.9887,143.1027) --
      (532.9887,146.9893) -- (533.4460,153.6193) -- (536.4181,156.5914) --
      (539.8474,156.5914) -- (541.6764,159.7922) -- (545.1057,160.2494) --
      (548.9923,165.9650) -- (556.0796,170.0802) -- (558.1372,172.8236) --
      (559.0517,180.2539) -- (559.7376,183.5689) -- (562.0238,185.1693) --
      (562.2524,186.5410) -- (560.1948,189.9703) -- (560.4234,193.1711) --
      (562.9383,197.0576) -- (565.4531,198.2007) -- (568.4252,198.6580) --
      (569.7676,200.0381) -- (615.0659,197.3687) -- cycle;

    \path[USA map/state, USA map/NC, local bounding box=NC] (834.9815,294.3155) -- (837.0665,299.2329) -- (840.6231,305.6993) --
      (843.0480,308.1242) -- (843.6946,310.3875) -- (841.2697,310.5491) --
      (842.0780,311.1958) -- (841.7547,315.3989) -- (839.1681,316.6922) --
      (838.5215,318.7938) -- (837.2282,321.7037) -- (833.5100,323.3203) --
      (831.0851,322.9970) -- (829.6301,322.8353) -- (828.0135,321.5420) --
      (828.3369,322.8353) -- (828.3369,323.8053) -- (830.2768,323.8053) --
      (831.0851,325.0986) -- (829.1452,331.4033) -- (833.3483,331.4033) --
      (833.9950,333.0199) -- (836.2582,330.7567) -- (837.5515,330.2717) --
      (835.6116,333.8282) -- (832.5400,338.6781) -- (831.2468,338.6781) --
      (830.1151,338.1931) -- (827.3669,338.8397) -- (822.1938,341.2646) --
      (815.7273,346.5994) -- (812.3325,351.2876) -- (810.3926,357.7540) --
      (809.9076,360.1789) -- (805.2194,360.6639) -- (799.7663,362.0005) --
      (789.8199,353.7980) -- (777.2103,346.2000) -- (774.3004,345.3916) --
      (761.6909,346.8466) -- (757.4145,347.5967) -- (755.7979,344.3635) --
      (752.8275,342.2468) -- (736.3381,342.7318) -- (729.0634,343.5401) --
      (720.0104,348.0666) -- (713.8673,350.6532) -- (692.6897,353.2398) --
      (693.1898,349.1854) -- (694.9681,347.7305) -- (697.7163,347.0838) --
      (698.3630,343.3656) -- (702.5661,340.6174) -- (706.4460,339.1624) --
      (710.6492,335.6059) -- (715.0140,333.5043) -- (715.6606,330.4328) --
      (719.5405,326.5529) -- (720.1871,326.3913) .. controls (720.1871,326.3913) and
      (720.1871,327.5229) .. (720.9955,327.5229) .. controls (721.8038,327.5229) and
      (722.9354,327.8462) .. (722.9354,327.8462) -- (725.1986,324.2897) --
      (727.3002,323.6430) -- (729.5635,323.9664) -- (731.1801,320.4098) --
      (734.0900,317.8232) -- (734.5750,315.7217) -- (734.7625,312.0735) --
      (739.0390,312.0510) -- (746.2375,311.1952) -- (761.9948,308.9427) --
      (777.1308,306.8562) -- (798.7713,302.1368) -- (818.7546,297.8782) --
      (829.9316,295.4724) -- (834.9815,294.3156) -- cycle(839.2520,327.5221) --
      (841.8386,325.0164) -- (844.9909,322.4298) -- (846.5267,321.7831) --
      (846.6884,319.7624) -- (846.0417,313.6193) -- (844.5868,311.2752) --
      (843.9401,309.4161) -- (844.6676,309.1736) -- (847.4159,314.6701) --
      (847.8200,319.1157) -- (847.6584,322.5106) -- (844.2635,324.0464) --
      (841.4344,326.4713) -- (840.3028,327.6838) -- (839.2520,327.5221) -- cycle;

    \path[USA map/state, USA map/DC, local bounding box=DC] (805.8194,250.8438) -- (803.9612,249.0197) -- (802.7285,248.3334) --
      (804.1715,246.3109) -- (807.0606,248.2594) -- (805.8194,250.8438) -- cycle;

    \path[USA map/state, USA map/MA, local bounding box=MA] (899.6235,173.2539) -- (901.7954,172.5681) -- (902.2527,170.8534) --
      (903.2815,170.9677) -- (904.3103,173.2539) -- (903.0529,173.7112) --
      (899.1662,173.8255) -- (899.6235,173.2539) -- cycle(890.2499,174.0541) --
      (892.5362,171.4250) -- (894.1365,171.4250) -- (895.9655,172.9110) --
      (893.5650,173.9398) -- (891.3931,174.9686) -- (890.2499,174.0541) --
      cycle(855.4508,152.0659) -- (873.0977,147.4253) -- (875.3609,146.7786) --
      (877.2750,143.9829) -- (881.0118,142.3196) -- (883.9010,146.7324) --
      (881.4761,151.9056) -- (881.1528,153.3605) -- (883.0927,155.9471) --
      (884.2244,155.1388) -- (886.0026,155.1388) -- (888.2659,157.7253) --
      (892.1457,163.7068) -- (895.7023,164.1918) -- (897.9655,163.2218) --
      (899.7438,161.4435) -- (898.9355,158.6953) -- (896.8339,157.0787) --
      (895.3789,157.8870) -- (894.4090,156.5937) -- (894.8939,156.1087) --
      (896.9955,155.9471) -- (898.7738,156.7554) -- (900.7137,159.1803) --
      (901.6837,162.0902) -- (902.0070,164.5151) -- (897.8038,165.9700) --
      (893.9240,167.9099) -- (890.0441,172.4364) -- (888.1042,173.8914) --
      (888.1042,172.9214) -- (890.5291,171.4665) -- (891.0141,169.6882) --
      (890.2058,166.6167) -- (887.2959,168.0716) -- (886.4876,169.5266) --
      (886.9726,171.7898) -- (884.9063,172.7902) -- (882.1591,168.2631) --
      (878.7642,163.8983) -- (876.6937,162.0858) -- (870.1604,163.9620) --
      (865.0681,165.0128) -- (844.3929,169.6050) -- (843.7252,164.8371) --
      (844.3718,154.2484) -- (848.6611,153.3592) -- (855.4508,152.0659) -- cycle;

    \path[USA map/state, USA map/TN, local bounding box=TN] (696.6779,318.2541) -- (644.7848,323.2656) -- (629.0252,325.0439) --
      (624.4040,325.5566) -- (620.5357,325.5289) -- (620.3147,329.6297) --
      (612.1293,329.8937) -- (605.1779,330.5403) -- (597.0871,330.4165) --
      (595.6733,337.4894) -- (593.9771,342.9694) -- (590.6839,345.7202) --
      (589.3352,350.1013) -- (589.0118,352.6879) -- (584.9703,354.9511) --
      (586.4253,358.5076) -- (585.4553,362.8725) -- (584.4869,363.6621) --
      (692.6455,353.2546) -- (693.0487,349.2996) -- (694.8595,347.8093) --
      (697.6936,347.0598) -- (698.3656,343.3428) -- (702.4642,340.6379) --
      (706.5111,339.1438) -- (710.5947,335.5735) -- (715.0308,333.5480) --
      (715.5520,330.4807) -- (719.6166,326.4957) -- (720.1674,326.3815) .. controls
      (720.1674,326.3815) and (720.1986,327.5132) .. (721.0070,327.5132) .. controls
      (721.8153,327.5132) and (722.9469,327.8677) .. (722.9469,327.8677) --
      (725.2101,324.2799) -- (727.2805,323.6333) -- (729.5556,323.9285) --
      (731.1539,320.3956) -- (734.1092,317.7517) -- (734.5308,315.8126) --
      (734.8398,312.1015) -- (732.6932,311.9017) -- (730.0916,313.9300) --
      (723.0983,313.9591) -- (704.7390,316.3460) -- (696.6779,318.2542) -- cycle;

    \path[USA map/state, USA map/AR, local bounding box=AR] (593.8248,343.0530) -- (589.8449,343.7697) -- (584.7327,343.1356) --
      (585.1534,341.5336) -- (588.1332,338.9669) -- (589.0766,335.3106) --
      (587.2476,332.3385) -- (508.8300,334.8534) -- (510.4304,341.7121) --
      (510.4304,349.9425) -- (511.8021,360.9165) -- (512.0307,398.7534) --
      (514.3170,400.6967) -- (517.2891,399.3250) -- (520.0325,400.4681) --
      (520.7129,407.0414) -- (576.3341,405.9008) -- (577.4798,403.8104) --
      (577.1932,400.2609) -- (575.3675,397.2888) -- (576.9662,395.8036) --
      (575.3675,393.2921) -- (576.0517,390.7822) -- (577.4201,385.1768) --
      (579.9383,383.1142) -- (579.2524,380.8296) -- (582.9104,375.4578) --
      (585.6539,374.0894) -- (585.5404,372.5959) -- (585.1949,370.7702) --
      (588.0519,365.1715) -- (590.4549,363.9149) -- (590.8391,360.4873) --
      (592.6097,359.2456) -- (589.4662,358.7613) -- (588.1248,354.7509) --
      (590.9288,352.3742) -- (591.4791,350.3550) -- (592.7586,346.3083) --
      (593.8248,343.0530) -- cycle;

    \path[USA map/state, USA map/MO, local bounding box=MO] (558.4402,248.1132) -- (555.9203,245.0259) -- (554.7772,242.7397) --
      (490.4200,245.1402) -- (488.1337,245.2545) -- (489.3912,247.7694) --
      (489.1626,250.0556) -- (491.6774,253.9422) -- (494.7638,258.0574) --
      (497.8502,260.8009) -- (500.0114,261.0295) -- (501.5082,261.9440) --
      (501.5082,264.9161) -- (499.6792,266.5164) -- (499.2219,268.8027) --
      (501.2795,272.2320) -- (503.7944,275.2041) -- (506.3092,277.0331) --
      (507.6810,288.6928) -- (507.9951,324.7650) -- (508.2237,329.4518) --
      (508.6810,334.8353) -- (531.1140,333.9685) -- (554.3200,333.2826) --
      (575.1246,332.4816) -- (586.7794,332.2513) -- (588.9488,335.6773) --
      (588.2646,338.9848) -- (585.1773,341.3878) -- (584.6050,343.2252) --
      (589.9834,343.6824) -- (593.8784,342.9966) -- (595.5956,337.5029) --
      (596.2470,331.6461) -- (598.3450,329.0910) -- (600.9411,327.6041) --
      (600.9925,324.5538) -- (602.0085,322.6174) -- (600.3143,320.0736) --
      (598.9833,321.0579) -- (596.9907,318.8306) -- (595.7057,314.0716) --
      (596.5067,311.5534) -- (594.5626,308.1258) -- (592.7319,303.5500) --
      (587.9325,302.7506) -- (580.9637,297.1519) -- (579.2449,293.0383) --
      (580.0442,289.8376) -- (582.1035,283.7799) -- (582.5624,280.9163) --
      (580.6133,279.8850) -- (573.7579,279.0873) -- (572.7299,277.3752) --
      (572.6181,273.1448) -- (567.1312,269.7138) -- (560.1557,261.9423) --
      (557.8695,254.6264) -- (557.6392,250.4011) -- (558.4402,248.1132) -- cycle;

    \path[USA map/state, USA map/GA, local bounding box=GA] (672.2923,355.5518) -- (672.2923,357.7342) -- (672.4539,359.8358) --
      (673.1006,363.2307) -- (676.4955,371.1521) -- (678.9204,381.0134) --
      (680.3753,387.1565) -- (681.9919,392.0063) -- (683.4469,398.9577) --
      (685.5485,405.2625) -- (688.1350,408.6574) -- (688.6200,412.0522) --
      (690.5599,412.8605) -- (690.7216,414.9621) -- (688.9433,419.8119) --
      (688.4584,423.0452) -- (688.2967,424.9851) -- (689.9133,429.3499) --
      (690.2366,434.6847) -- (689.4283,437.1096) -- (690.0750,437.9179) --
      (691.5299,438.7262) -- (691.7346,441.9443) -- (693.9676,445.2939) --
      (696.2181,447.4559) -- (704.1395,447.6176) -- (714.9592,446.9709) --
      (736.4716,445.6777) -- (741.9173,445.0033) -- (746.4946,445.0310) --
      (746.6562,447.9409) -- (749.2428,448.7492) -- (749.5661,444.3843) --
      (747.9495,439.8578) -- (749.0811,438.2412) -- (754.9009,439.0495) --
      (759.8783,439.3673) -- (759.1029,433.0685) -- (761.3661,423.0456) --
      (762.8211,418.8424) -- (762.3361,416.2558) -- (765.6705,410.0115) --
      (765.1602,408.6599) -- (763.2468,409.3644) -- (760.6602,408.0711) --
      (760.0136,405.9695) -- (758.7203,402.4130) -- (756.4571,400.3114) --
      (753.8705,399.6648) -- (752.2539,394.8150) -- (749.3289,388.4800) --
      (745.1257,386.5400) -- (743.0241,384.6001) -- (741.7308,382.0135) --
      (739.6292,380.0736) -- (737.3660,378.7803) -- (735.1027,375.8704) --
      (732.0312,373.6072) -- (727.5047,371.8289) -- (727.0197,370.3740) --
      (724.5948,367.4641) -- (724.1098,366.0091) -- (720.7149,361.0386) --
      (717.1951,361.1378) -- (713.4401,358.7817) -- (712.0219,357.4884) --
      (711.6985,355.7102) -- (712.5693,353.7702) -- (714.7960,352.6601) --
      (714.1620,350.5629) -- (672.2923,355.5518) -- cycle;

    \path[USA map/state, USA map/SC, local bounding box=SC] (764.9433,408.1649) -- (763.1662,409.1344) -- (760.5796,407.8411) --
      (759.9330,405.7395) -- (758.6397,402.1830) -- (756.3765,400.0814) --
      (753.7899,399.4347) -- (752.1733,394.5849) -- (749.4251,388.6035) --
      (745.2219,386.6635) -- (743.1203,384.7236) -- (741.8270,382.1370) --
      (739.7254,380.1971) -- (737.4622,378.9038) -- (735.1989,375.9939) --
      (732.1274,373.7307) -- (727.6009,371.9524) -- (727.1159,370.4975) --
      (724.6910,367.5876) -- (724.2060,366.1326) -- (720.8111,360.9595) --
      (717.4162,361.1211) -- (713.3747,358.6962) -- (712.0814,357.4029) --
      (711.7581,355.6247) -- (712.5664,353.6848) -- (714.8297,352.7148) --
      (714.3189,350.4257) -- (720.0870,348.0891) -- (729.2025,343.5001) --
      (736.9772,342.6918) -- (753.0916,342.2693) -- (755.7298,344.1468) --
      (757.4089,347.5050) -- (761.7113,346.8950) -- (774.3208,345.4400) --
      (777.2307,346.2484) -- (789.8402,353.8464) -- (799.9483,361.9681) --
      (794.5272,367.4264) -- (791.9406,373.5695) -- (791.4556,379.8743) --
      (789.8390,380.6826) -- (788.7074,383.4308) -- (786.2825,384.0775) --
      (784.1809,387.6340) -- (781.4327,390.3822) -- (779.1694,393.7771) --
      (777.5528,394.5854) -- (773.9963,397.9803) -- (771.0864,398.1419) --
      (772.0564,401.3751) -- (767.0449,406.8716) -- (764.9433,408.1649) -- cycle;

    \path[USA map/state, USA map/KY, local bounding box=KY] (725.9944,295.2707) -- (723.7011,297.6724) -- (720.1229,301.6664) --
      (715.1983,307.1311) -- (713.9826,308.8469) -- (713.9201,310.9484) --
      (709.5402,313.1125) -- (703.8821,316.5074) -- (696.6502,318.3063) --
      (644.7823,323.2051) -- (629.0228,324.9834) -- (624.4016,325.4961) --
      (620.5332,325.4684) -- (620.3063,329.6887) -- (612.1269,329.8332) --
      (605.1755,330.4799) -- (597.1880,330.4197) -- (598.3958,329.0996) --
      (600.8953,327.5587) -- (601.1239,324.3580) -- (602.0384,322.5290) --
      (600.4316,319.9901) -- (601.2334,318.0833) -- (603.4967,316.3051) --
      (605.5983,315.6584) -- (608.3465,316.9517) -- (611.9030,318.2450) --
      (613.0347,317.9217) -- (613.1963,315.6584) -- (611.9030,313.2335) --
      (612.2264,310.9702) -- (614.1663,309.5153) -- (616.7529,308.8687) --
      (618.3695,308.2220) -- (617.5612,306.4437) -- (616.9145,304.5038) --
      (618.4211,303.5080) .. controls (618.4241,303.4709) and (619.6751,299.9857) ..
      (619.6594,299.8502) -- (622.7127,298.3715) -- (628.0324,297.4016) --
      (632.5265,296.9166) -- (633.9189,298.5440) -- (635.4472,299.4148) --
      (637.0380,296.3066) -- (640.2250,295.0240) -- (642.4301,296.5080) --
      (642.8407,297.5071) -- (644.0142,297.2430) -- (643.8525,294.2901) --
      (646.9834,292.5409) -- (649.1315,291.4674) -- (650.6609,293.1283) --
      (653.9790,293.0841) -- (654.5663,291.5128) -- (654.1988,289.2496) --
      (656.7994,285.2511) -- (661.5759,281.8132) -- (662.2819,276.9773) --
      (665.2069,276.5214) -- (668.9983,274.8757) -- (671.4417,273.1675) --
      (671.2433,271.6025) -- (670.1009,270.1476) -- (670.6667,267.1527) --
      (674.8516,267.0352) -- (677.1515,266.2894) -- (680.4989,267.7185) --
      (682.5530,272.0833) -- (687.6853,272.0941) -- (689.7363,274.3023) --
      (691.3517,274.1546) -- (693.9534,272.8765) -- (699.1905,273.4498) --
      (701.7654,273.6673) -- (703.4530,271.6111) -- (706.0709,270.1852) --
      (707.9527,269.4781) -- (708.5993,272.3147) -- (710.6428,273.3731) --
      (713.2855,275.4556) -- (713.4030,281.1288) -- (714.2113,282.7012) --
      (716.8010,284.2575) -- (717.5727,286.5520) -- (721.7325,289.9890) --
      (723.5379,293.6122) -- (725.9944,295.2707) -- cycle;

    \path[USA map/state, USA map/AL, local bounding box=AL] (631.3065,460.4157) -- (629.8159,446.0942) -- (627.0676,427.3416) --
      (627.2293,413.2771) -- (628.0376,382.2382) -- (627.8759,365.5872) --
      (628.0410,359.1681) -- (672.5255,355.5487) -- (672.3777,357.7311) --
      (672.5394,359.8327) -- (673.1860,363.2276) -- (676.5809,371.1489) --
      (679.0058,381.0102) -- (680.4607,387.1534) -- (682.0773,392.0032) --
      (683.5323,398.9546) -- (685.6339,405.2593) -- (688.2205,408.6542) --
      (688.7054,412.0491) -- (690.6454,412.8574) -- (690.8070,414.9590) --
      (689.0287,419.8088) -- (688.5438,423.0420) -- (688.3821,424.9820) --
      (689.9987,429.3468) -- (690.3220,434.6816) -- (689.5137,437.1065) --
      (690.1604,437.9148) -- (691.6153,438.7231) -- (691.9435,441.6119) --
      (686.3458,441.2584) -- (679.5561,441.9050) -- (654.0137,444.8149) --
      (643.6021,446.2217) -- (643.3807,449.0991) -- (645.1590,450.8774) --
      (647.7456,452.8173) -- (648.3264,460.7527) -- (642.7844,463.3256) --
      (640.0361,463.0023) -- (642.7844,461.0624) -- (642.7844,460.0924) --
      (639.7128,454.1110) -- (637.4496,453.4643) -- (635.9946,457.8291) --
      (634.7013,460.5774) -- (634.0547,460.4157) -- (631.3065,460.4157) -- cycle;

    \path[USA map/state, USA map/LA, local bounding box=LA] (607.9671,459.1612) -- (604.6824,455.9951) -- (605.6924,450.4949) --
      (605.0310,449.6018) -- (595.7693,450.6084) -- (570.7410,451.0673) --
      (570.0568,448.6726) -- (570.9696,440.2169) -- (574.2855,434.2711) --
      (579.3169,425.5800) -- (578.7428,423.1820) -- (579.9994,422.5012) --
      (580.4583,420.5487) -- (578.1721,418.4927) -- (578.0603,416.5503) --
      (576.2296,412.2048) -- (576.0826,405.8662) -- (520.6088,406.7902) --
      (520.6374,416.3637) -- (521.3233,425.7373) -- (522.0091,429.6238) --
      (524.5240,433.7390) -- (525.4385,438.7688) -- (529.7823,444.2557) --
      (530.0109,447.4564) -- (530.6968,448.1423) -- (530.0109,456.6013) --
      (527.0388,461.6310) -- (528.6392,463.6886) -- (527.9533,466.2035) --
      (527.2675,473.5194) -- (525.8957,476.7201) -- (526.0182,480.3365) --
      (530.7047,478.8164) -- (542.8180,479.0234) -- (553.1643,482.5799) --
      (559.6307,483.7116) -- (563.3489,482.2566) -- (566.5821,483.3882) --
      (569.8153,484.3582) -- (570.6236,482.2566) -- (567.3904,481.1250) --
      (564.8038,481.6100) -- (562.0556,479.9934) .. controls (562.0556,479.9934) and
      (562.2173,478.7001) .. (562.8639,478.5384) .. controls (563.5105,478.3768) and
      (565.9355,477.5685) .. (565.9355,477.5685) -- (567.7137,479.0234) --
      (569.4920,478.0534) -- (572.7252,478.7001) -- (574.1801,481.1250) --
      (574.5035,483.3882) -- (579.0299,483.7116) -- (580.8082,485.4898) --
      (579.9999,487.1064) -- (578.7066,487.9147) -- (580.3232,489.5313) --
      (588.7296,493.0879) -- (592.2861,491.7946) -- (593.2561,489.3697) --
      (595.8426,488.7230) -- (597.6209,487.2681) -- (598.9142,488.2381) --
      (599.7225,491.1479) -- (597.4592,491.9562) -- (598.1059,492.6029) --
      (601.5008,491.3096) -- (603.7640,487.9147) -- (604.5723,487.4298) --
      (602.4707,487.1064) -- (603.2790,485.4898) -- (603.1174,484.0349) --
      (605.2189,483.5499) -- (606.3506,482.2566) -- (606.9972,483.0649) .. controls
      (606.9972,483.0649) and (606.8355,486.1365) .. (607.6439,486.1365) .. controls
      (608.4522,486.1365) and (611.8470,486.7831) .. (611.8470,486.7831) --
      (615.8885,488.7230) -- (616.8585,490.1780) -- (619.7684,490.1780) --
      (620.9000,491.1479) -- (623.1633,488.0764) -- (623.1633,486.6214) --
      (621.8700,486.6214) -- (618.4751,483.8732) -- (612.6553,483.0649) --
      (609.4221,480.8017) -- (610.5537,478.0534) -- (612.8170,478.3768) --
      (612.9786,477.7301) -- (611.2004,476.7602) -- (611.2004,476.2752) --
      (614.4336,476.2752) -- (616.2119,473.2036) -- (614.9186,471.2637) --
      (614.5953,468.5155) -- (613.1403,468.6771) -- (611.2004,470.7787) --
      (610.5537,473.3653) -- (607.4822,472.7186) -- (606.5122,470.9404) --
      (608.2905,469.0005) -- (610.1938,465.5548) -- (609.1327,463.1426) --
      (607.9671,459.1612) -- cycle;

    \path[USA map/state, USA map/MS, local bounding box=MS] (631.5588,459.3446) -- (631.3046,460.6007) -- (626.1314,460.6007) --
      (624.6765,459.7924) -- (622.5749,459.4691) -- (615.7851,461.4090) --
      (614.0069,460.6007) -- (611.4203,464.8039) -- (610.3178,465.5819) --
      (609.1939,463.0939) -- (608.0508,459.2074) -- (604.6215,456.0066) --
      (605.7646,450.4621) -- (605.0787,449.5476) -- (603.2498,449.7762) --
      (595.3318,450.6496) -- (570.7853,451.0230) -- (570.0156,448.7976) --
      (570.8890,440.4208) -- (574.0058,434.7480) -- (579.2329,425.6031) --
      (578.7871,423.1705) -- (580.0240,422.5142) -- (580.4599,420.5948) --
      (578.1424,418.5158) -- (578.0273,416.3743) -- (576.1915,412.2532) --
      (576.0825,406.2905) -- (577.4101,403.8095) -- (577.1868,400.3937) --
      (575.4173,397.3111) -- (576.9437,395.8289) -- (575.3731,393.3294) --
      (575.8303,391.6772) -- (577.4077,385.1508) -- (579.8937,383.1145) --
      (579.2520,380.7475) -- (582.9100,375.4450) -- (585.7419,374.0885) --
      (585.5209,372.4134) -- (585.2328,370.7323) -- (588.1088,365.1646) --
      (590.4545,363.9331) -- (590.6062,363.0401) -- (627.9496,359.1589) --
      (628.1345,365.4422) -- (628.2962,382.0933) -- (627.4879,413.1322) --
      (627.3262,427.1966) -- (630.0744,445.9493) -- (631.5588,459.3446) -- cycle;

    \path[USA map/state, USA map/IA, local bounding box=IA] (569.1915,199.5843) -- (569.4559,202.3705) -- (571.6796,202.9478) --
      (572.6336,204.1731) -- (573.1336,206.0285) -- (576.9264,209.3871) --
      (577.6123,211.7786) -- (576.9380,215.2031) -- (575.3556,218.4351) --
      (574.5563,221.1768) -- (572.3836,222.7789) -- (570.6680,223.3513) --
      (565.0890,225.2115) -- (563.6976,229.0602) -- (564.4262,230.4319) --
      (566.2667,232.1145) -- (565.9838,236.1508) -- (564.2206,237.6887) --
      (563.4492,239.3318) -- (563.5764,242.1081) -- (561.6901,242.5654) --
      (560.0647,243.6703) -- (559.7859,245.0229) -- (560.0647,247.1378) --
      (558.5137,248.2539) -- (556.0431,245.1206) -- (554.7806,242.6707) --
      (489.0447,245.1856) -- (488.1267,245.3510) -- (486.0743,240.8351) --
      (485.8457,234.2050) -- (484.2453,230.0898) -- (483.5595,224.8315) --
      (481.2732,221.1735) -- (480.3588,216.3724) -- (477.6153,208.8279) --
      (476.4722,203.4552) -- (475.1004,201.2833) -- (473.5001,198.5399) --
      (475.4541,193.6960) -- (476.8258,187.9805) -- (474.0823,185.9229) --
      (473.6251,183.1794) -- (474.5396,180.6645) -- (476.2542,180.6645) --
      (558.9082,179.3951) -- (559.7425,183.5782) -- (561.9947,185.1392) --
      (562.2514,186.5622) -- (560.2219,189.9516) -- (560.4123,193.1571) --
      (562.9271,196.9553) -- (565.4539,198.2489) -- (568.5332,198.7519) --
      (569.1915,199.5843) -- cycle;

    \path[USA map/state, USA map/MN, local bounding box=MN] (475.2378,128.8244) -- (474.7806,120.3653) -- (472.9516,113.0494) --
      (471.1226,99.5607) -- (470.6654,89.7299) -- (468.8364,86.3006) --
      (467.2360,81.2709) -- (467.2360,70.9829) -- (467.9219,67.0963) --
      (466.1009,61.6446) -- (496.2334,61.6799) -- (496.5567,53.4352) --
      (497.2033,53.2735) -- (499.4666,53.7585) -- (501.4065,54.5668) --
      (502.2148,60.0633) -- (503.6697,66.2064) -- (505.2863,67.8230) --
      (510.1362,67.8230) -- (510.4595,69.2779) -- (516.7642,69.6012) --
      (516.7642,71.7028) -- (521.6141,71.7028) -- (521.9374,70.4095) --
      (523.0690,69.2779) -- (525.3322,68.6313) -- (526.6255,69.6012) --
      (529.5354,69.6012) -- (533.4153,72.1878) -- (538.7501,74.6127) --
      (541.1750,75.0977) -- (541.6599,74.1277) -- (543.1149,73.6428) --
      (543.5999,76.5526) -- (546.1864,77.8459) -- (546.6714,77.3609) --
      (547.9647,77.5226) -- (547.9647,79.6242) -- (550.5513,80.5942) --
      (553.6228,80.5942) -- (555.2394,79.7859) -- (558.4726,76.5526) --
      (561.0592,76.0677) -- (561.8675,77.8459) -- (562.3525,79.1392) --
      (563.3224,79.1392) -- (564.2924,78.3309) -- (573.1837,78.0076) --
      (574.9620,81.0791) -- (575.6086,81.0791) -- (576.3223,79.9949) --
      (580.7622,79.6242) -- (580.1501,81.9037) -- (576.2113,83.7408) --
      (566.9656,87.8019) -- (562.1908,89.8088) -- (559.1193,92.3954) --
      (556.6944,95.9519) -- (554.4311,99.8318) -- (552.6529,100.6401) --
      (548.1264,105.6515) -- (546.8331,105.8132) -- (542.5053,108.5703) --
      (540.0424,111.7754) -- (539.8138,114.9668) -- (539.9082,123.0102) --
      (538.5322,124.6989) -- (533.4506,128.4589) -- (531.2205,134.4413) --
      (534.0923,136.6750) -- (534.7722,139.9020) -- (532.9169,143.1409) --
      (533.0877,146.8889) -- (533.4566,153.6193) -- (536.4848,156.6213) --
      (539.8138,156.6213) -- (541.7050,159.7539) -- (545.0841,160.2572) --
      (548.9433,165.9287) -- (556.0306,170.0454) -- (558.1737,172.9205) --
      (558.8449,179.3600) -- (477.6334,180.5048) -- (477.2955,144.8280) --
      (476.8382,141.8559) -- (472.7230,138.4265) -- (471.5799,136.5976) --
      (471.5799,134.9972) -- (473.6375,133.3968) -- (475.0092,132.0251) --
      (475.2379,128.8244) -- cycle;

    \end{scope}
}

\tikzstyle{grapharrow}=[->,line width=1.0pt, shorten >=1pt]
\tikzset{
  textnode/.style={
	minimum height=15pt,inner sep=0pt, draw,circle, line width=1.0pt
  },
}

\usepackage{tabularx,environ}
\makeatletter
\newcommand{\problemtitle}[1]{\gdef\@problemtitle{#1}}
\newcommand{\probleminput}[1]{\gdef\@probleminput{#1}}
\newcommand{\problemquestion}[1]{\gdef\@problemquestion{#1}}
\NewEnviron{feasproblem}{
  \problemtitle{}\probleminput{}\problemquestion{}
  \BODY
  \par\addvspace{.5\baselineskip}
  \noindent
  \begin{tabularx}{\textwidth}{@{\hspace{\parindent}} l X c}
    \multicolumn{2}{@{\hspace{\parindent}}l}{\@problemtitle} \\
    \textbf{Input:} & \@probleminput \\
    \textbf{Output:} & \@problemquestion
  \end{tabularx}
  \par\addvspace{.5\baselineskip}
}
\makeatother

\usepackage[figuresleft]{rotating}
\usepackage{caption}

\usepackage{multirow}

\usepackage[T1]{fontenc}
\usepackage{listings}

\usepackage{enumitem}
\newlist{paraenum}{enumerate*}{1}
\setlist[paraenum]{label=(\emph{\roman*})}

\usepackage{tipa} 

\usepackage{bm}

\usepackage{pifont}

\usepackage{todonotes}
\usepackage[flushleft]{threeparttable}
\usepackage{multicol}

\usepackage[english]{babel}
\usepackage{blindtext}

\begin{document}

\begin{frontmatter}

\title{The Traveling Tournament Problem: An Overview}

\author[add1,add2]{David Van Bulck}
\ead{david.vanbulck@ugent.be}

\author[add3]{Fan Yang\corref{mycorrespondingauthor}}
\cortext[mycorrespondingauthor]{Corresponding author}
\ead{fan_yang@shnu.edu.cn}

\author[add1,add2]{Dries Goossens}
\ead{dries.goossens@ugent.be}

\author[add4]{Michael Trick}
\ead{trick@cmu.edu}

\address[add1]{Faculty of Economics and Business Administration, Ghent University, Ghent, Belgium}
\address[add2]{FlandersMake@UGent -- core lab CVAMO, Ghent, Belgium}
\address[add3]{School of Finance and Business, Shanghai Normal University, Shanghai, China}
\address[add4]{Carnegie Mellon University in Qatar}

\begin{center}
\emph{This paper is dedicated to Dr.\ Kelly Easton who was instrumental in defining the Traveling Tournament Problem, in providing early results on the problem, and in using those results in real-world sports schedules. Her tragically early passing in February 2026 was a tremendous loss both personally and for the sports scheduling field.}
\end{center}

\begin{abstract}
Over the past 25 years, the Traveling Tournament Problem (TTP) has become one of the most extensively studied optimization problems in sports scheduling. At its core, the TTP seeks to minimize the total travel distance incurred by teams that travel directly between opponents' venues during consecutive away games. The problem originated from the scheduling challenges faced by Major League Baseball, where it was identified as the central computational difficulty. This paper provides a comprehensive overview of the literature on the TTP. We review the principal problem variants and benchmark instances, and summarize the current state of the art in lower bounds, approximation guarantees, and exact and heuristic optimization algorithms. Moreover, we contribute to the continued development of the field by tracking and validating lower and upper bounds, while succeeding the repository originally established by Prof.\ Michael Trick as part of the RobinX sports timetabling project.
Finally, we identify several open questions and outline promising directions for future research.
\end{abstract}


\end{frontmatter}
\section{Introduction}

In 1995, a seemingly simple question emerged from Major League Baseball (MLB): could operations research techniques be used to create a timetable for its league? 
What followed was a decade-long struggle that did not result in a schedule for MLB (that would come later), but did give rise to one of the most extensively studied optimization problems in sports scheduling: the Traveling Tournament Problem (TTP). 
The TTP, formally introduced in the seminal paper of \citet{Easton2001}, asks for a compact double round-robin tournament that minimizes the overall distance traveled by the teams.
That is, over the course of the season, each of an even number of teams plays every other team once at home and once away, and it plays exactly one game per round.
In order to minimize travel, teams are allowed to group away games into a road trips: when consecutively playing away, teams move directly between the opponent's venues, without returning home.
Two additional constraints make the problem realistic. 
First, to avoid excessive travel fatigue, long periods without home games for the fans, and too many consecutive home games, no team is allowed to play more than a given number of consecutive home or away games.
Second, for reasons of game attractiveness and fairness, no two teams may play each other in consecutive rounds.
For an example of the input of the problem, a round-robin timetable solution, and the individual road trips of a particular team, we refer to \Cref{fig:exampleTTP}.

Although the problem is deceptively simple to describe, it is notoriously difficult to solve. 
To understand part of its difficulty, observe that each team in itself solves a vehicle routing problem: using its own venue as the `depot', it visits every other team in road trips with a given maximal length (the `capacity' of the vehicles), while minimizing travel distance. 
Yet the teams individual routing problems cannot be solved separately, because every team plays a game in every round, so the trips of one team dictate which other teams have to be home and thus cannot be on a road trip themselves.
And as schedulers discovered the hard way, timetables are inherently fragile: even a minor modification in one game triggers a cascade of changes that may ultimately compromise the quality or even the feasibility of the resulting schedule. 
As a result, the largest TTP real-distance instance solved to optimality contains merely 10 teams.
This is a surprisingly small number, especially in light of the vehicle routing problem where instances involving hundreds of customers are now routinely solved to optimality.

\begin{figure}
\centering
\begin{subfigure}[b]{0.48\textwidth}
\centering
\begin{tikzpicture}[scale=0.35]
\tiny
\tikzset{textnode/.style={minimum height=10pt, inner sep=0pt, draw, circle, line width=0.7pt}}

\USA[every state={draw=white, ultra thick, fill=black!10}]

\node[textnode] (MOa) at (NY)[yshift=18pt,xshift=9pt]{\textsc{mo}};
\node[textnode] (NYa) at (NY)[yshift=-6pt,xshift=10pt]{\textsc{ny}};
\node[textnode] (PHa) at (PA)[xshift=15pt,yshift=-5pt]{\textsc{ph}};
\node[textnode] (ATa) at (GA)[xshift=-4pt,yshift=7pt]{\textsc{at}};
\node[textnode] (FLa) at (FL)[xshift=13pt,yshift=1pt]{\textsc{fl}};
\node[textnode] (PIa) at (PA)[xshift=-6pt,yshift=0pt]{\textsc{pi}};
\node[textnode] (CIa) at (OH)[xshift=-6pt,yshift=-4pt]{\textsc{ci}};
\node[textnode] (CHa) at (IL)[xshift=1pt,yshift=10pt]{\textsc{ch}};
\node[textnode] (STa) at (IL)[xshift=-1pt,yshift=-4pt]{\textsc{st}};
\node[textnode] (MIa) at (WI)[xshift=3pt,yshift=-4pt]{\textsc{mi}};

\end{tikzpicture}
\caption{Input. Pairwise travel distances between the teams, here visualized by the distances between the teams' home locations}
\label{fig:inputTTP}
\end{subfigure}
\hfill
\begin{subfigure}[b]{0.48\textwidth}
\centering
\addtocounter{subfigure}{1} 
\begin{tikzpicture}[scale=0.35]
\tiny
\tikzset{textnode/.style={minimum height=10pt, inner sep=0pt, draw, circle, line width=0.7pt}}

\USA[every state={draw=white, ultra thick, fill=black!10}]

\node[textnode] (MOc) at (NY)[yshift=18pt,xshift=9pt]{\textsc{mo}};
\node[textnode] (NYc) at (NY)[yshift=-6pt,xshift=10pt]{\textsc{ny}};
\node[textnode] (PHc) at (PA)[xshift=15pt,yshift=-5pt]{\textsc{ph}};
\node[textnode] (ATc) at (GA)[xshift=-4pt,yshift=7pt]{\textsc{at}};
\node[textnode] (FLc) at (FL)[xshift=13pt,yshift=1pt]{\textsc{fl}};
\node[textnode] (PIc) at (PA)[xshift=-6pt,yshift=0pt]{\textsc{pi}};
\node[textnode] (CIc) at (OH)[xshift=-6pt,yshift=-4pt]{\textsc{ci}};
\node[textnode] (CHc) at (IL)[xshift=1pt,yshift=10pt]{\textsc{ch}};
\node[textnode] (STc) at (IL)[xshift=-1pt,yshift=-4pt]{\textsc{st}};
\node[textnode] (MIc) at (WI)[xshift=3pt,yshift=-4pt]{\textsc{mi}};

\draw[->,line width=0.8pt, shorten >=2pt] (CIc) -- (PHc);
\draw[->,line width=0.8pt, shorten >=0.2pt] (PHc) -- (NYc);
\draw[->,line width=0.8pt, shorten >=2pt] (NYc) -- (MOc);
\path (MOc) edge[bend left=-30, line width=0.8pt,->, shorten >=2pt] (CIc);

\draw[<->,line width=0.8pt, shorten >=2pt, shorten <=2pt] (CIc) -- (PIc);

\draw[->,line width=0.8pt, shorten >=2pt] (CIc) -- (STc);
\draw[->,line width=0.8pt, shorten >=2pt] (STc) -- (FLc);
\draw[->,line width=0.8pt, shorten >=2pt] (FLc) -- (ATc);
\draw[->,line width=0.8pt, shorten >=2pt] (ATc) -- (CIc);

\draw[->,line width=0.8pt, shorten >=2pt] (CIc) -- (CHc);
\draw[->,line width=0.8pt, shorten >=0.2pt] (CHc) -- (MIc);
\draw[->,line width=0.8pt, shorten >=2pt] (MIc) -- (CIc);
\end{tikzpicture}
\caption{Visualization of the road trips for team CI in the timetable provided below}
\label{fig:tripsTTP}
\end{subfigure}

\medskip
\begin{subfigure}{\textwidth}
\addtocounter{subfigure}{-2 } 
\centering
\caption{An optimal timetable. The first team in parentheses is the home team, the second the away team. The away games of team CI are typeset in bold. The maximal sequence of home and away games is limited to three.}
\label{fig:solTTP}
\scriptsize
\setlength{\tabcolsep}{3pt}
\resizebox{\textwidth}{!}{%
\begin{tabular}{l*{18}{c}}
\toprule
 & R1 & R2 & R3 & R4 & R5 & R6 & R7 & R8 & R9 & R10 & R11 & R12 & R13 & R14 & R15 & R16 & R17 & R18\\
\midrule
1 & (AT,CH) & (AT,MI) & (AT,FL) & (NY,ST) & (NY,MO) & \textbf{(NY,CI)} & (AT,PH) & (AT,NY) & (AT,PI) & (NY,AT) & (NY,CH) & (NY,MI) & \textbf{(AT,CI)} & (AT,ST) & (AT,MO) & (NY,FL) & (NY,PI) & (NY,PH)\\
2 & (MO,PH) & (MO,ST) & (PH,ST) & (PH,MO) & \textbf{(PH,CI)} & (FL,PH) & \textbf{(MO,CI)} & (FL,PI) & (PH,NY) & (PH,CH) & (PH,MI) & (MO,CH) & (MO,MI) & (FL,MO) & (PH,FL) & (PH,PI) & (PH,AT) & (MO,PI)\\
3 & (FL,MI) & (FL,CH) & (MO,NY) & (CI,FL) & (CH,PI) & (PI,MO) & (FL,NY) & (CI,PH) & (CI,ST) & (CI,MO) & (MO,AT) & \textbf{(FL,CI)} & (FL,ST) & (CI,PI) & (PI,ST) & (CI,AT) & (MO,FL) & (FL,AT)\\
4 & \textbf{(PI,CI)} & (PI,PH) & (PI,CH) & (CH,AT) & (ST,FL) & (ST,AT) & (PI,MI) & (CH,MO) & (CH,FL) & (ST,PI) & (PI,FL) & (PI,AT) & (PI,NY) & (CH,NY) & (CI,CH) & (CH,MI) & \textbf{(CH,CI)} & (ST,CH)\\
5 & (ST,NY) & (CI,NY) & (CI,MI) & (MI,PI) & (MI,AT) & (MI,CH) & (CH,ST) & (MI,ST) & (MI,MO) & (MI,FL) & \textbf{(ST,CI)} & (ST,PH) & (CH,PH) & (MI,PH) & (MI,NY) & (ST,MO) & (ST,MI) & \textbf{(MI,CI)}\\
\bottomrule
\end{tabular}%
}
\end{subfigure}
\caption{Example of a problem instance with 10 teams (NL10) together with an optimal solution and a visualization of the road trips for one particular team}
\label{fig:exampleTTP}
\end{figure}

\begin{figure}
\centering
\footnotesize
\begin{tikzpicture}
\begin{axis}[
    ybar,
    ylabel={},
    xlabel={},
    bar width=25pt,
    xtick=data,
    xticklabels={2000-2004,2005-2009,2010-2014,2015-2019,2020-2024, 2025-},
    height=5cm,
    enlarge x limits=0.2,
    ymin=0,
    ymax=320,
    width=0.8\linewidth,
    axis line style={draw=none}, 
    ytick={},
    yticklabels={},
    major tick length=0pt 
]

\addplot[fill=blue!50, nodes near coords] coordinates {
    (2002.5, 51)
    (2007.5, 184)
    (2012.5, 281)
    (2017.5, 223)
    (2022.5, 258)
    (2027.5, 81)
};
\end{axis}
\end{tikzpicture}
\caption{Number of Google Scholar hits for `"Traveling Tournament Problem" OR "Travelling Tournament Problem"', per five-year periods}
\label{fig:citations}
\end{figure}
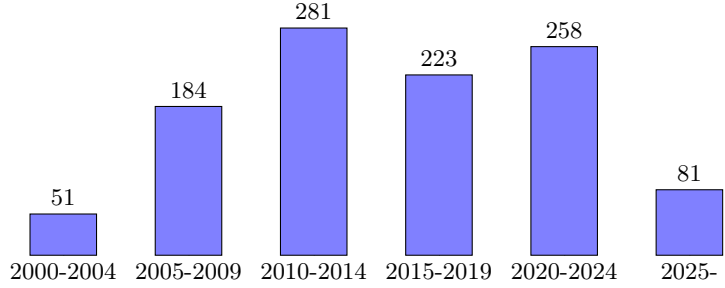

The TTP continues to attract a steady stream of research: over the past 25 years, it has accumulated more than a thousand Google Scholar hits (\Cref{fig:citations} gives the evolution per five-year period). 
In our view, the literature on the TTP can be organized into three broad research waves. 
During a first wave, roughly spanning the first decade of research (2000-2009), the focus was mostly on metaheuristics. 
This wave resulted in the first dedicated non-trivial neighborhoods for sports scheduling, respecting the elementary round-robin structure of the timetable and applicable more broadly than for the TTP alone.
Successful applications include simulated annealing, GRASP, and tabu search (see e.g., \citet{Anagnostopoulos2006,Ribeiro2007,Gaspero2007}).
In a second wave, roughly spanning 2008 to 2012, the emphasis shifted towards exact methods. 
Building on earlier branch-and-price ideas by \citet{Easton2003}, in 2008 \citet{Irnich2010} succeeded to prove optimality for the eight-team NL8 instance.
One year later, a provenly optimal solution to the NL10 instance was found by \citet{Uthus2009,Uthus2012} employing dept-first and best-first search algorithms.
A third and ongoing wave, starting around 2010, concerns the theoretical analysis of the problem, including the first complexity proofs for the unconstrained variant by \citet{Bhattacharyya2009,Bhattacharyya2016} and a growing body of results on approximation algorithms (e.g., \citet{Miyashiro2012,Westphal2012,Zhao2025e}) and inapproximability (e.g., \citet{Bendayan2023,Zhao2025}).

With this growing body of literature, we believe an overview of the field is timely.
Indeed, to the best of our knowledge, no comprehensive overview on the TTP exists. 
The aim of this paper is to fill this gap and to provide a guide for researchers who wish to start working on the TTP. 
The remainder of this paper is organized as follows. \Cref{sec:origins} recounts how the TTP came into existence. 
\Cref{sec:variants} introduces the formal problem definition, reviews the principal problem variants, and presents the benchmark instances that are commonly used in the literature. 
\Cref{sec:lowerBound} surveys the lower bounds that have been proposed for the TTP, followed by an overview of the approximation results in \Cref{sec:approx}. 
\Cref{sec:optimization} gives an overview of the exact and heuristic optimization algorithms that have been developed. 
Finally, \Cref{sec:conclusion} concludes and discusses directions for future research.
Throughout the paper, we identify seven open problems, drawn either from the literature or from questions that we believe deserve attention in the years to come.
Beyond the survey itself, this paper contributes to the community in a concrete way: it brings together all commonly used TTP instances and their best-known solutions, including previously lost instance classes as well as full solutions rather than merely objective values, and computes lower bounds for instances for which none or only weaker ones existed. Moreover, it takes over the maintenance of the repository started by Prof.\ Michael Trick, which now officially lives on as part of the RobinX sports timetabling project (see \url{www.robinxval.ugent.be/RobinX/travelRepo.php}).

\section{Origins of the traveling tournament problem}
\label{sec:origins}

Although some earlier papers address travel minimization for teams undertaking road trips inspired by real-world sports tournaments, e.g., \citet{Campbell1976} (basketball) and \citet{Russell1994} (baseball), the TTP is first formally defined by \citet{Easton2001}. The authors were inspired by a request to schedule the Major League Baseball (MLB), which Prof.\ Michael Trick received from Doug Bureman, an executive with a local MLB team, in 1995. In Trick's own words, the Traveling Tournament Problem originated as follows.

\begin{quote}
When Doug approached me in 1995, I was full of confidence that I could schedule Major League Baseball. 
I figured I could put together a quick greedy algorithm followed by local search, perhaps with simulated annealing (then about 10 years old) that would clearly get good answers reasonably quickly.  I knew that Henry and Holly Stephenson, baseball’s schedulers, used computers to count aspects of the schedule, but the bulk of the work was based on hand-generated patterns that they laboriously put together.  Thus, I was certain that my vast array of operations research techniques could create better schedules, faster.
But it was not the case. Very quickly, I realized that it was not straightforward to create a local search algorithm.  Given a schedule, simply exchanging the opponents on two games led to a cascading set of necessary changes that might not result in a feasible schedule.  And even creating a greedy heuristic seemed difficult: when scheduling week-by-week, it didn’t take long to reach a point where the schedule could not be completed.
Stymied on the heuristic front, I returned to an integer programming approach.  And here I received another shock:  the integer programs I formulated not only could not be solved to optimality, but the software never even came back with feasible integer solutions.  
When I created the integer programs, I included everything Major League Baseball required in a schedule. In addition to scheduling the 2,430 games, there were many, many additional requirements, including that every team must play half their weekends at home, that stadiums are sometimes unavailable for home games, and that teams do not want to travel too far.

I did learn some key aspects that made things a bit simpler.  In baseball, teams typically played 2 to 4 games against a single opponent (a series) before moving on to the next opponent.  So instead of 2,430 games, I was really trying to schedule 780 series.  Better, initially MLB had 2 leagues that did not play each other, so I was really trying to schedule a 14 team league and a 16 team league (with weak interaction due to the two team cities), not a 30 team league. 
I also learned that I had to have variables that represented ``trips'' that teams took, rather than individual series.  So if team A went to team B, then C, then D before returning home, then that road trip could be represented by a single integer variable.  
But, despite all I had learned, and continuing improvements in the underlying solver, I still got nowhere in solving the integer programs.  The solver would branch and branch and branch without improving the lower bound (where the objective was primarily to minimize distance traveled) and without finding a feasible solution for an upper bound.
Of course, being an academic, this failure to find anything approaching a useful approach didn’t stop me from presenting the work at various conferences and workshops.  More importantly, I began working with George Nemhauser, a famous person in integer programming circles, and Kelly Easton (then a doctoral student) on scheduling other leagues.  It turned out that while MLB was out of reach, integer programming methods could work wonders for smaller leagues, and we began scheduling college sports leagues.  The schedules we created were much better than the hand-schedules being used, and quickly multiple college leagues played our schedules.
But Major League Baseball continued to haunt us.  In 1999, I did what I should have done much earlier: I started to explore why my integer programs were hard to solve.  I decided to start removing constraints until I got down to a solvable problem.  I took away the weekend constraints, and then the team requests, and so on and so on.  Surprisingly, I realized I could take away almost everything and still end up with an unsolvable problem.  In fact, even when I reduced the problem to a double round-robin (instead of the quadruple round-robin that then made up the MLB schedule), as long as I tried to find a schedule with minimum distance traveled, I could not even solve a six team instance to optimality!   The problem was not with all the “complicating” constraints: the difficulty was with the underlying basic problem of finding minimum distance round-robins.
So the Traveling Tournament Problem was born: how to create such minimum distance round robins.  
We first formally presented this work as a short paper at the Constraint Programming Conference in 2001 \citep{Easton2001}. I was very worried that someone would quickly find a way to solve this problem, making the TTP irrelevant.  As history would show, despite dozens of papers, no one has found a quick way to solve the instances we presented, and even finding provably optimal solutions to 12 team round robins remains beyond current capabilities.
Twenty-five years after the TTP was first developed, it remains a rich source of inspiration for research in heuristic and optimization approaches.
For Major League Baseball, ten years after my initial forays into the problem, our group finally created schedules that MLB would play.  We would do so for almost all the next 12 years, before being supplanted by firms that continued to develop optimization approaches for this problem.
\end{quote}

\section{The TTP, problem variants, and instances}
\label{sec:variants}

\subsection{Preliminaries}

Let $T$ be the set of teams, where $|T|=n\geq 4$ is even, and $R$ the set of rounds (often corresponding to weekends; also called time slots).
During each round, a team can play at most one game.
Furthermore, we define with $D$ an $n\times n$ symmetric distance matrix, where element $d_{i,j}\in D$ reflects the non-negative integer distance between the home venue of team $i$ and team $j$.
It is common in the TTP literature on lower bounds and approximation ratio's to assume that $D$ respects the triangle inequality (i.e., $d_{i,j} + d_{j,k} \geq d_{i,k}\, \forall i,j,k\in T$), and that $d_{i,i}=0$ for all $i\in T$.
Nonetheless, there are benchmark instances that violate this assumption.
A game is an ordered pair of teams $(i,j)$ in which $i\in T$ is the home team providing the venue where the game is played, and $j\in T\setminus \{i\}$ is the away team.
If a team plays two or more away games in a row, it has a road trip and travels directly from the venue of one opponent to the venue of the next opponent.
We define a leg of a road trip as a single direct journey between two venues.
For instance, the road trip of team CI to teams CH and MI in \Cref{fig:tripsTTP} consists of legs CI-CH, CH-MI, and MI-CI.
Similarly, a team has a home stand when it plays two or more consecutive games at home.
It is assumed that all teams are initially at home, and that teams need to return home after playing their last away game.
A double round-robin tournament (2RR) is a collection of games in which every team plays every other team exactly once at home and once away, and is called compact if it uses the minimum number of rounds needed (i.e., $2n-2$ since we assume $n$ even). Finally, the complete graph on $n$ vertices is denoted by $K_n$.

\subsection{Classic version of the problem}
\label{subsec:classic}

The classic version of the TTP is defined as follows.
{
\begin{feasproblem}
	\problemtitle{\textbf{\textsc{Traveling Tournament  with Maximal Trip Length $k$ (TTP($k$))}}}
	\probleminput{Teams $T$ and rounds $R$ with $|R|=2n-2$, distance matrix $D$, and a non-negative integer $k$.}
	\problemquestion{
		A compact 2RR timetable minimizing the sum of the distances travelled by all teams such that:
		\begin{description}
			\item[C1] The length of home stands and road trips is at most $k$, and 
			\item[C2] No game $(i,j)$ is immediately followed by game $(j,i)$.
		\end{description}
	}
\end{feasproblem}
}

Constraints C1 and C2 are referred to as the `at-most' and `no-repeater' constraints, respectively.
Unless otherwise specified, this paper assumes that $k=3$, which is by far the most studied case in the literature.
Since no compact timetable exists that perfectly alternates between home and away games for all teams (see e.g., \citet{Werra1981}), it is clear that $k$ must be greater than 1.
The version of the problem without the at-most constraint (i.e., $k=n-1$) is referred to as the unconstrained traveling tournament problem (UTTP).
We note that the original description of the TTP by \citet{Easton2001} also includes a parameter to regulate the minimal length of home stands and road trips.
However, this constraint has been ignored in the literature so far.

With regard to the complexity status of TTP($k$), the following results are known (for an overview of results, see \Cref{tab:complexity}).
The first complexity proof is provided by \citet{Bhattacharyya2009,Bhattacharyya2016}, who show that the UTTP is $\mathcal{NP}$-hard by providing a reduction from the (1,2)-Traveling Salesman Problem (TSP).
\citet{Thielen2011} provide a reduction from 3-SAT to TTP(3), showing that the TTP remains $\mathcal{NP}$-hard for $k=3$.
\citet{Chatterjee2021} generalizes this result for any fixed $k>3$.
The complexity status of TTP(2), though, is still open.

\begin{open}
    Determine the computational complexity status of TTP(2), i.e., TTP(2)$\in \mathcal{P}$?
\end{open}

On the positive side, several approximation algorithms are known for TTP(2) with an approximation ratio which goes to 1 as $n$ increases to infinity (the earliest of which is due to \citet{Thielen2012}). 
%
%
This implies that TTP(2) has a PTAS (see \citet{Zhao2025}).
In contrast, \citet{Bendayan2023} show that the UTTP is APX-hard, thus no PTAS exists unless $\mathcal{P}=\mathcal{NP}$, solving an open problem raised by \citet{Imahori2014}.
This result has later been generalized to arbitrary $k>2$ by \citet{Zhao2025}.

\begin{table}
	\scriptsize
	\centering
	\begin{tabular}{l llll}
		\toprule
		Publications  			    & TTP(2)  & TTP(3) 	& TTP($k>3$) 	& UTTP ($k=n-1$)\\
		\midrule \addlinespace[5pt]
        \textbf{P vs.\ NP-hard}\\
                                    & ?\\
		\citet{Thielen2011}     	&& NP-hard\\
		\citet{Chatterjee2021}     	&&& NP-hard\tnote{\tiny 1}\\
		\citet{Bhattacharyya2016}   &&&&NP-hard\\[5pt]
        \textbf{PTAS vs.\ APX-hard}\\
        \citet{Thielen2012} & PTAS\\
        \citet{Zhao2025} && APX-hard & APX-hard\\
        \citet{Bendayan2023}    	&&&&APX-hard\\
		\bottomrule
	\end{tabular}
	\caption{Overview of complexity results.}
	\label{tab:complexity}
\end{table}

\subsection{Problem variants}

There are several variants of the TTP that consider alternative constraints and tournament formats. 
One variant that has been particularly researched is the mirrored TTP (m-TTP), introduced by \citet{Ribeiro2007}.
In this variant, the 2RR is divided into two phases, where the second phase is identical to the first except that the home and away teams are swapped for every match.
The bipartite traveling tournament problem extends the TTP to competitions consisting of two leagues with dedicated and consecutive time slots for interleague games, where the interleague play is represented by a bipartite tournament (see \citet{Hoshino2011b}). 
\citet{Hoshino2011d} propose the multi-round TTP which generalizes the TTP to $l$-RR tournaments where teams meet each other exactly $l>2$ times (in the TTP, $l=2$). 
Long before the introduction of the TTP, \citet{Russell1994} already consider a variant that can now be seen as the TTP(2), however, focusing on a 1RR (i.e., $l=1$) and additionally enforcing that each team plays approximately half of its games at home.
\citet{Devriesere2026} introduce the TTP in the context of incomplete round-robin tournaments, in which fewer than $n-1$ time slots are available and, consequently, not all teams can play each other; they term this variant the incomplete TTP.
In contrast, \citet{Bao2010} assume that there are more than $2n-2$ time slots, meaning that teams no longer play on every time slot.
This variant is known as the time-relaxed TTP (RTTP).
Several other variants of the TTP assume that part of the timetable has already been fixed. 
For instance, the TTP with predefined-venues (TTP-PV) requires constructing a travel-minimal 1RR schedule for the second half of a phased tournament, where each pair of teams must play at the other team's venue compared to their first match (see \citet{Melo2009}). 
The timetable constrained distance TTP seeks an optimal home-away assignment when the opponents of each team in each time slot are given (see \citet{Rasmussen2008}). 
The TTP with trip preferences (TP-TTP) aims to create a travel-minimal timetable while ensuring that all road trips belong to a predefined set of trips.
A real-life application of the TP-TTP in the context of Argentina's professional basketball league is discussed in \citet{Duran2019}.
\citet{Osicka2023} introduce fairness objectives to the TTP: rather than solely minimizing the overall distance traveled, they use a cooperative game theory approach to fairly distribute the travel distance over the teams.

Finally, there is a large strand of literature that focuses on the Traveling Umpire Problem (TUP), which aims to assign umpires to games, given a timetable for the tournament. The goal is to minimize umpire travel while ensuring that no umpire handles the game of a particular team too frequently (see \citet{Trick2012}). 
\citet{Bender2016} propose an integration of the TTP with the TUP.

\subsection{Common problem instances}

The popularity of the TTP is partly attributable to the availability of several well-known benchmark instances (see \Cref{tab:TTP_Instances}).
Perhaps most known are the National League (NLx) problem instances proposed by \citet{Easton2001}, which are based on the actual locations of Major League Baseball teams in the National League.
Other instance classes that are based on real-life tournaments include the National Football League (NFLx, see \citet{Uthus2009}), Super 14 Rugby League (SUPx, see \citet{Uthus2009}), and Brazilian soccer championship (BRAx, \citet{Ribeiro2007}). In addition, the Galaxy instances (GALx, see \citet{Uthus2012}) embed the team locations in a 3D-coordinate space, where the distance matrix is based on the number of light-years between stars in the universe.

\begin{table}
	\centering
	\footnotesize
	\caption{Overview of existing traveling tournament instance classes}
	\label{tab:TTP_Instances}
	\begin{tabular}{r l c c}
		\toprule
			Abbreviation & Full name & Type & No.\@ teams\\
			\midrule
			CONx 	& Constant distance instances				& Artificial 	& 	[4,40]\\
			CIRCx & Circular distance instances				& Artificial 	& 	[4,40]\\
			LINEx & Linear distance instances 				& Artificial 	&   [4,40]\\
			INCRx & Increasing distance instances 			& Artificial 	&   [4,40]\\
			GALx 	& Galaxy instances							& Artificial 	& 	[4,40]\\
			NLx 	& National League instances 		 		& Real-life 	& 	[4,16]\\
			NFLx 	& National Football League instances 		& Real-life 	& 	[4,32]\\
			SUPx 	& Super 14 Rugby League instances 	 		& Real-life 	& 	[4,14]\\
			BRAx 	& Brazilian soccer championship instances 	& Real-life 	& 	[4,24]\\
		\bottomrule
	\end{tabular}
\end{table}

\begin{figure}
\footnotesize
\begin{subfigure}[b]{0.24\textwidth}
\centering
\begin{tikzpicture}[scale=0.7]
\def\radius{1.5cm}

\foreach \i/\label in {1/1, 2/2, 3/3, 4/4}
\node[circle, draw, fill=white] (\label) at ({-90*(\i+2)}:\radius) {\label};

\draw (1) edge[bend left=28] node[above]{1} (2);
\draw (2) edge[bend left=28] node[right]{1} (3);
\draw (3) edge[bend left=28] node[below]{1} (4);
\draw (4) edge[bend left=28] node[left]{1} (1);
\end{tikzpicture}
\caption{CIRC4}
\end{subfigure}
\hfill
\begin{subfigure}[b]{0.48\textwidth}
\centering
\begin{tikzpicture}[scale=0.5]
  \node[circle, draw, fill=white] (n1) at (0,3) {1};
  \node[circle, draw, fill=white] (n2) at (2,3) {2};
  \node[circle, draw, fill=white] (n3) at (4,3) {3};
  \node[circle, draw, fill=white] (n4) at (6,3) {4};

  \draw (n1) -- (n2) node[midway, above]{1};
  \draw (n2) -- (n3) node[midway, above]{1};
  \draw (n3) -- (n4) node[midway, above]{1};

  \node[circle, draw, fill=white] (n1bis) at (0,0) {1};
  \node[circle, draw, fill=white] (n2bis) at (2,0) {2};
  \node[circle, draw, fill=white] (n3bis) at (6,0) {3};
  \node[circle, draw, fill=white] (n4bis) at (12,0) {4};

  \draw (n1bis) -- (n2bis) node[midway, above]{1};
  \draw (n2bis) -- (n3bis) node[midway, above]{2};
  \draw (n3bis) -- (n4bis) node[midway, above]{3};
\end{tikzpicture}
\caption{LINE4 and INCR4}
\end{subfigure}
\hfill
\begin{subfigure}[b]{0.24\textwidth}
\centering
\begin{tikzpicture}
    
  \node[circle, draw, fill=white] (v1) at (0,0) {1};
  \node[circle, draw, fill=white] (v2) at (2,0) {2};
  \node[circle, draw, fill=white] (v3) at (1,{sqrt(3)}) {3};
  \node[circle, draw, fill=white] (v4) at (1,{sqrt(3)/2-0.18}) {4};
  
  \draw (v1) -- (v3) node[midway, left] {1};
  \draw (v1) -- (v4) node[midway, right, xshift=0pt, yshift=-2pt] {1};
  \draw (v2) -- (v4) node[midway, left, xshift=0pt, yshift=-2pt] {1};
  \draw (v1) -- (v2) node[midway, below] {1};
  \draw (v2) -- (v3) node[midway, right] {1};
  \draw (v3) -- (v4) node[midway, left, xshift=2pt] {1};
\end{tikzpicture}
\caption{CON4}
\end{subfigure}
\caption{Graph representation of artificial TTP problem instances for which the associated TSP is trivial. Distance $d_{i,j}$ is equal to the shortest path length between node $i$ and $j$ in the graph.}
\label{fig:graph}
\end{figure}
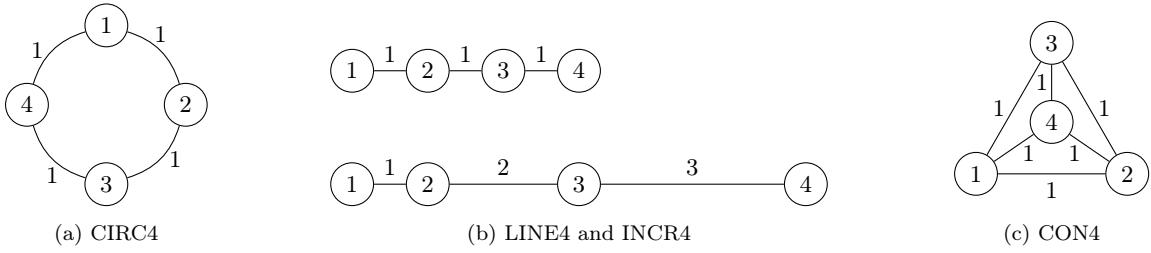

The above instance classes pose significant challenges: the largest instance solved to proven optimality contains only ten teams.
What if we add some topological structure to the instances, such that the associated TSP becomes trivial?
In this light, \citet{Easton2001} introduce the circular distance instances (CIRCx) where teams are positioned on a circle and travel to other teams by moving along it (see also \Cref{fig:graph}).
Following a similar idea, \citet{Hoshino2012} propose linear (LINEx) and incremental (INCRx) distance instances, where all teams are located on a straight line.
In LINEx the distance between adjacent teams is always 1, whereas in INCRx the distance between adjacent teams increases by one unit at each step.
However, even with this simplified structure, none of the instances with more than 10 teams has been solved to optimality so far.
What if we completely ignore the distance matrix, by having a distance of 1 between any pair of teams?
The resulting instances are known as the constant distance instances (CONx, see \citet{Urrutia2006}), and all optimal solutions are known for up to 16 teams.
We refer to the TTP with constant distances as the CD-TTP.

For more than twenty years, the website maintained by Prof.\ Michael Trick\footnote{See \url{https://mat.tepper.cmu.edu/TOURN/}.} has tracked the best-known solutions for several of the above problem instances. 
As part of this project, this website has been transferred to RobinX\footnote{See \url{www.robinxval.ugent.be/RobinX/}.} (see \citet{VanBulck2019}).
RobinX offers many advantages over a plain-text only website, including XML-based storage of the problem instances and solutions (rather than just their objective values), and an online validator for checking solution feasibility.
Moreover, missing instances have been included and several solutions and lower bounds have been restored or even improved.
For instance, we have included the LINEx and INCRx instance classes and in collaboration with \citeauthor{Melo2009} we have retrieved several TTP-PV problem instances that had been lost and were therefore not available on Trick's website. 
Moreover, lower bounds have been computed for the first time for several larger problem instances that had previously been considered computationally too demanding (see the next section).

\section{Lower Bounds}
\label{sec:lowerBound}

\begin{figure}
  \centering
  
   \begin{subfigure}[t]{0.3\textwidth}
    \centering
\begin{tikzpicture}[scale=0.5]
\footnotesize



\node[textnode] (MO3) at (NY)[yshift=26pt,xshift=8pt]{\textsc{mo}};
\node[textnode] (NY3) at (NY)[yshift=-8pt,xshift=14pt]{\textsc{ny}};
\node[textnode] (PH3) at (PA)[xshift=22pt,yshift=-7.5pt]{\textsc{ph}};
\node[textnode] (AT3) at (GA)[xshift=-5pt,yshift=10pt]{\textsc{at}};
\node[textnode] (FL3) at (FL)[xshift=18pt,yshift=2pt]{\textsc{fl}};
\node[textnode] (PI3) at (PA)[xshift=-8pt,yshift=2pt]{\textsc{pi}};
\node[textnode] (CI3) at (OH)[xshift=-8pt,yshift=-6pt]{\textsc{ci}};
\node[textnode] (CH3) at (IL)[xshift=2pt,yshift=14pt]{\textsc{ch}};
\node[textnode] (ST3) at (IL)[xshift=-2pt,yshift=-5pt]{\textsc{st}};
\node[textnode] (MI3) at (WI)[xshift=4pt,yshift=-5pt]{\textsc{mi}};

\draw[line width=1.2pt] (FL3) -- (AT3);
\draw[line width=1.2pt] (PH3) -- (NY3);
\draw[line width=1.2pt] (CI3) -- (ST3);
\draw[line width=1.2pt] (CH3) -- (MI3);
\draw[line width=1.2pt] (PI3) -- (MO3);

\draw[line width=1.2pt, dotted] (CI3) -- (MO3);
\path (CI3) edge[bend left=-16, line width=1.2pt,dotted] (NY3);
\draw[line width=1.2pt, dotted] (CI3) -- (PH3);
\draw[line width=1.2pt, dotted] (CI3) -- (AT3);
\path (CI3) edge[bend left=16, line width=1.2pt,dotted] (FL3);
\draw[line width=1.2pt, dotted] (CI3) -- (PI3);
\draw[line width=1.2pt, dotted] (CI3) -- (CH3);
\draw[line width=1.2pt, dotted] (CI3) -- (ST3);
\draw[line width=1.2pt, dotted] (CI3) -- (MI3);
\end{tikzpicture}
    \caption{ILB(2) = matching}
    \label{fig:ilb2}
  \end{subfigure}
  \hfill
  \begin{subfigure}[t]{0.3\textwidth}
    \centering
    
\begin{tikzpicture}[scale=0.5]
\footnotesize



\node[textnode] (MO3) at (NY)[yshift=26pt,xshift=8pt]{\textsc{mo}};
\node[textnode] (NY3) at (NY)[yshift=-8pt,xshift=14pt]{\textsc{ny}};
\node[textnode] (PH3) at (PA)[xshift=22pt,yshift=-7.5pt]{\textsc{ph}};
\node[textnode] (AT3) at (GA)[xshift=-5pt,yshift=10pt]{\textsc{at}};
\node[textnode] (FL3) at (FL)[xshift=18pt,yshift=2pt]{\textsc{fl}};
\node[textnode] (PI3) at (PA)[xshift=-8pt,yshift=2pt]{\textsc{pi}};
\node[textnode] (CI3) at (OH)[xshift=-8pt,yshift=-6pt]{\textsc{ci}};
\node[textnode] (CH3) at (IL)[xshift=2pt,yshift=14pt]{\textsc{ch}};
\node[textnode] (ST3) at (IL)[xshift=-2pt,yshift=-5pt]{\textsc{st}};
\node[textnode] (MI3) at (WI)[xshift=4pt,yshift=-5pt]{\textsc{mi}};

\draw[line width=1.2pt] (CI3) -- (MO3);
\draw[line width=1.2pt] (NY3) -- (MO3);
\draw[line width=1.2pt] (NY3) -- (PH3);
\draw[line width=1.2pt] (CI3) -- (PH3);
33
\draw[line width=1.2pt] (CI3) -- (AT3);
\draw[line width=1.2pt] (AT3) -- (FL3);
\draw[line width=1.2pt] (FL3) -- (PI3);
\draw[line width=1.2pt] (PI3) -- (CI3);
33
\draw[line width=1.2pt] (CH3) -- (MI3);
\draw[line width=1.2pt] (MI3) -- (CI3);
\draw[line width=1.2pt] (CH3) -- (ST3);
\draw[line width=1.2pt] (CI3) -- (ST3);


\end{tikzpicture}
    \caption{ILB(3) = VRP}
    \label{fig:ilb3}
  \end{subfigure}
  \hfill
  \begin{subfigure}[t]{0.3\textwidth}
    \centering
\begin{tikzpicture}[scale=0.3]
\footnotesize



\node[textnode] (MO3) at (NY)[yshift=26pt,xshift=13pt]{\textsc{mo}};
\node[textnode] (NY3) at (NY)[yshift=-8pt,xshift=14pt]{\textsc{ny}};
\node[textnode] (PH3) at (PA)[xshift=22pt,yshift=-7.5pt]{\textsc{ph}};
\node[textnode] (AT3) at (GA)[xshift=-5pt,yshift=10pt]{\textsc{at}};
\node[textnode] (FL3) at (FL)[xshift=18pt,yshift=2pt]{\textsc{fl}};
\node[textnode] (PI3) at (PA)[xshift=-8pt,yshift=15pt]{\textsc{pi}};
\node[textnode] (CI3) at (OH)[xshift=-8pt,yshift=-6pt]{\textsc{ci}};
\node[textnode] (CH3) at (IL)[xshift=2pt,yshift=14pt]{\textsc{ch}};
\node[textnode] (ST3) at (IL)[xshift=-2pt,yshift=-5pt]{\textsc{st}};
\node[textnode] (MI3) at (WI)[xshift=4pt,yshift=-5pt]{\textsc{mi}};

\draw[line width=1.2pt] (CI3) -- (CH3);
\draw[line width=1.2pt] (CH3) -- (MI3);
\path (ST3) edge[bend left=30, line width=1.2pt] (MI3);
\draw[line width=1.2pt] (ST3) -- (AT3);
\draw[line width=1.2pt] (AT3) -- (FL3);
\draw[line width=1.2pt] (FL3) -- (PH3);
\draw[line width=1.2pt] (PH3) -- (NY3);
\draw[line width=1.2pt] (NY3) -- (MO3);
\path (MO3) edge[bend left=-30, line width=1.2pt] (PI3);
\draw[line width=1.2pt] (PI3) -- (CI3);


\end{tikzpicture}
    \caption{ILB($n-1$)=TSP}
    \label{fig:ilbn}
  \end{subfigure}
  \caption{Illustration of computation for the Independent Lower Bounds (ILB) for team CI}
  \label{fig:ilb}
\end{figure}
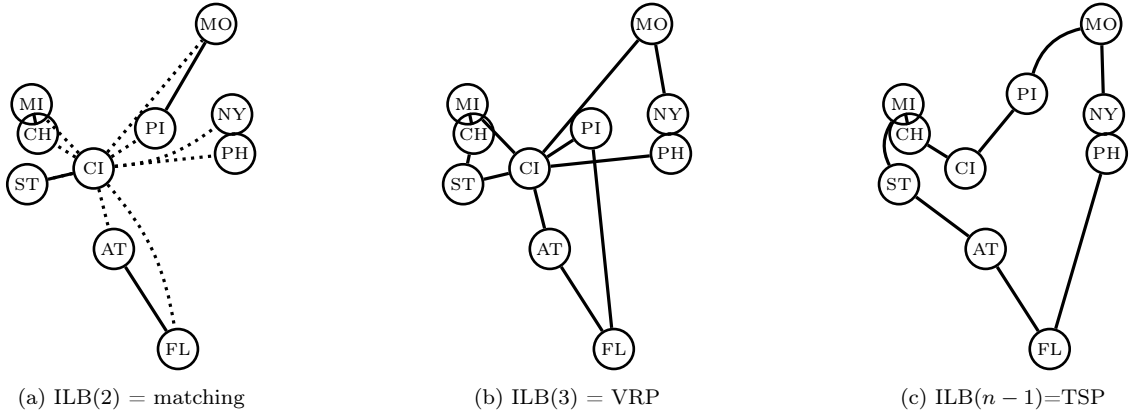

\paragraph{Independent lower bounds}

Perhaps the most commonly used technique to obtain lower bounds in TTP is to sum over the minimum distance that each team needs to travel to visit every other team, while taking into account the maximum length of road trips but ignoring all other scheduling constraints. 
This approach results in the `independent lower bound' (ILB), so named because it breaks down the TTP instance into $n$ sub-problems that can be solved independently. 

\citet{Campbell1976} and \citet{Ball1977} independently observe that a maximum trip length of two ($k=2$) requires every team to travel between their own venue and their opponent's venue at least once (see the dotted edges in \Cref{fig:ilb2}). 
In other words, summing over all teams, the total distance covered to move to the first opponent of each road trip of length one or two and to return home at the end of each road trip of length two is given by $\Delta=\sum_{i,j\in T:i\neq j} d_{i,j}$.
Disregarding the dotted edges, it becomes clear that a lower bound on the sum of travel from the first to the second opponent in each road trip of length two and returning home from a road trip of length one is given by a minimum cost matching, which can be found in polynomial time. 

\begin{proposition}[\textbf{ILB(2),} \citet{Campbell1976}]
\label{prop:ttp2}
Denote by $\gamma$ the minimum weight of a perfect matching in $K_n$.
Every solution of the TTP(2) has a total length of at least $\Delta+n\gamma$.
\end{proposition}

In contrast, for $k>2$, \citet{Easton2003b} proves that computing the ILB is $\mathcal{NP}$-complete in the strong sense for any constant $k$ by providing a reduction from the $\mathcal{NP}$-complete problem of partitioning a graph into isomorphic subgraphs restricted to paths.

As noted in \citet{Urrutia2007}, the ILB can be computed by solving $n$ independent vehicle routing problems (VRP).
More in particular, in the VRP of team $i\in T$, the venue of team $i$ serves as the depot, all other teams act as customers with a demand of one, and each vehicle has a capacity of $k$ (see \Cref{fig:ilb3}).

\begin{proposition}[\textbf{ILB($k$),} \citet{Easton2003}]
\label{prop:ilb}
Let $\pi_i$ be the optimal solution value for the VRP associated with team $i$.
Every solution of the TTP(k) has a total length of at least $\sum_{i\in T} \pi_i$.
\end{proposition}

Considering the UTTP, computing the ILB becomes equivalent to solving a TSP (see \Cref{fig:ilbn}).
This results in the following lower bound, which is obviously also valid for TTP(k) as the UTTP is a relaxation.

\begin{proposition}[\textbf{ILB($n-1$),} \citet{Yamaguchi2011}]
    \label{lb:tsp}
	Let $\rho$ be the length of an optimal TSP tour in $K_n$. 
	Every solution of the UTTP has a total length of at least $n\rho$.
\end{proposition}

An interesting conjecture by \citet{Bao2010} that, to the best of our knowledge, is still open is whether ILB(3) becomes tight when there are sufficiently more time slots than games per team.
\begin{open}
\label{open:rttp}
Determine whether there exists an $\alpha$ such that the ILB($k$) is tight for RTTP($k$) with $2n-2 + \alpha$ time slots, where $\alpha$ may depend on $n$ and $k$.
\end{open}

We have found the following results in the literature that are related to \Cref{open:rttp}.
Using at most two byes per team, \citet{Campbell1976} and \citet{Ball1977} show how to construct a schedule for RTTP(3) that attains ILB(2).
In contrast, computational experiments by \citet{Brandao2014} show that there exist instances for which two byes per team do not suffice for RTTP(3) to attain ILB(3).

Despite the fact that the independent lower bound can be quite strong in practice (for ILB(3) within 5\% of the optimal value, according to \cite{Miyashiro2012}), it is generally not tight. 
This is due to cases where the individual teams' road trips cannot be combined into a feasible timetable.
For instance, with regard to TTP(2), \citet{Thielen2012} prove that ILB(2) cannot be attained if the perfect matching used to compute the bound is unique.

\paragraph{Minimum number of legs}
A special type of problem instances for which stronger lower bounds are known is the constant distance instances of the CD-TTP.
In these instances, minimizing the overall travel distance is equivalent to minimizing the overall number of legs teams undertake during road trips. 
While \citet{Russell1994} already show a relation between travel minimization and break maximization (a team is said to have a break if it plays two consecutive home games or two consecutive away games), \citet{Urrutia2007} are the first to prove break maximization is equivalent to travel leg minimization.

\begin{proposition}[{\citet{Urrutia2006}}]
        Denote by $\beta$ and $\tau$ the total number of breaks and travel legs in a 2RR timetable, respectively.
        It holds that $\tau=n(2n-2)-\beta/2$.
\end{proposition}

As a consequence, any upper bound on the maximum number of breaks in a 2RR timetable can be transformed into a lower bound on the minimum number of legs in the timetable. 
Indeed, a team can have at most $2n-4$ breaks in a 2RR schedule, which occurs when the team plays all its home games before playing all its away games (or vice versa).
Moreover, since no two teams can play according to the same home-away pattern (see \citet{Werra1981}), the maximal number of breaks in a 2RR timetable is $2(2n-4)+(n-2)(2n-5)=2n^2-5n+2$, which gives rise to the following proposition.

\begin{proposition}[\citet{Urrutia2006}]
    Every solution of the TTP(k) contains at least $n^2+n/2-1$ legs.
\end{proposition}

\citet{Urrutia2006} provide a constructive proof to show that this bound is tight for the UTTP, implying that the UTTP restricted to constant distance instances is solvable in polynomial time.
Stronger bounds on the number of legs for TTP(3) are given by \citet{Fujiwara2007}.

\begin{proposition}[\citet{Fujiwara2007}]
Every solution of the TTP(3) has at least $(4/3)n^2-n$ legs if $n\equiv 0 \bmod{3}$, $(4/3)n^2-5/6n-1$ legs if $n\equiv 1 \bmod{3}$, and $(4/3)n^2-2/3n$ legs if $n\equiv 2 \bmod{3}$.
\end{proposition}

These bounds are tight for all instances of $n\equiv 1 \bmod{3}, n\leqslant 50$.
For similar lower bounds for the mirrored CD-TTP, we refer to \citet{Urrutia2006}. 
We note that \citet{Rasmussen2007} and \citet{VanBulck2023} further improve upon these bounds by considering logic-based and traditional Benders' cuts. 
Moreover, when minimizing travel legs, \citet{VanBulck2023} conjecture that it suffices to consider timetables where every team is paired with another that plays the exact opposite sequence of home and away games (referred to as a complementary Home/Away Pattern (HAP) set).

\begin{open}
Determine whether there always exists an optimal solution to CD-TTP(3) that has a complementary HAP set.
\end{open}

\paragraph{Minimum number of legs lower bound}

\citet{Urrutia2007} use optimal solution values (or lower bounds thereof) for the CD-TTP to strengthen the ILB.
In particular, they propose to simultaneously compute the ILB for all teams, with the additional constraint that the overall number of legs made by all teams is at least equal to the optimal value in the associated CD-TTP instance.
Even though its computation is considerably more challenging, the resulting bound, called the Minimum Number of Legs Lower Bound (MNLLB), often improves upon the ILB by several percentage points.

\begin{proposition}
	Let $\mu_n$ be a lower bound on the minimum number of legs in a 2RR timetable with $n$ teams, and let $\tau_i$ be the total number of legs made by team $i$.
	Every solution of the TTP(k) has a total length of at least $\min_{\sum_{i\in T}\tau_i \geq \mu_n} \big( \sum_{i\in T}\pi_i\big)$.
\end{proposition}

The MNLLB is further enhanced by \citet{Cheung2009}, who add the requirement that a set of logic-based Benders' cuts must be satisfied.
Nevertheless, this complicates the computation even further, limiting the author to compute bounds only for m-TTP(3).
During computational experiments with this code, however, we observe that the alternative IP formulation to compute the MNLLB provided by \citet{Cheung2009} is substantially faster than the IP model provided by \citet{Urrutia2007}.
Together with the author, this enabled us to compute the MNLLB for several problem instances for which it was not computed before, resulting in several new best lower bounds, all of which have been added to the RobinX website.

\paragraph{Other lower bounds}
Several other lower bounds have been proposed in the literature on approximation algorithms, which are typically easier to compute. 
Although these bounds often perform worse than the ILB in practice, evaluating the worst-case behavior of an algorithm based on these bounds can still provide useful insights. 
As an overview of these bounds is currently missing, we provide it here.

As each team needs to travel to the venue of every other team, a trivial lower bound on the total travel distance per team is the minimum spanning tree $\sigma$ in $K_n$. 
Recalling that $\gamma$ denotes the minimum weight of a perfect matching in $K_n$, \citet{Imahori2021} shows that $n\sigma+n\gamma$ is also a valid bound for TTP(2), but this bound is weaker than \Cref{prop:ttp2}.
Intuitively, spanning trees are cycle-free and thus ignore the fact that teams need to return home.
\citet{Yamaguchi2011} exploit this to provide the bound in \Cref{lb:mst} (for a variant, see \citet{Zhao2025e}).

\begin{proposition}[{\citet{Yamaguchi2011}}]
    \label{lb:mst}
	Let $\sigma$ be the length of a minimum spanning tree in $K_n$ and, given a solution, denote by $d_i^{\text{home}}$ the total distances covered by team $i$ to leave and return to its home before and after every road trip, respectively. Every solution of the TTP(k) has a total length of at least $n\sigma+\sum_{i}d_i^{\text{home}}/2$.
\end{proposition}

Alternatively, considering again the distances covered by the teams to leave and return to their home venues, \citet{Yamaguchi2011} propose the following bound (recall $\Delta$ denotes the sum of all pairwise distances).
\begin{proposition}[\citet{Yamaguchi2011}]
	Every solution of the TTP(k) has a total length of at least $\frac{2}{k-2}(\Delta-\sum_{i}d_i^{\text{home}})$.
\end{proposition}

\citet{Miyashiro2012} prove that the optimal value of TTP(3) is larger than or equal to $(2/3)\Delta$. 
\citet{Westphal2012} generalize this result for arbitrary $k$.

\begin{proposition}[\citet{Westphal2012}]
	\label{thm:maxDist}
	Every solution of the TTP(k) has a total length of at least $\frac{2 \Delta}{k}$.
\end{proposition}

For completeness, we note that several techniques have been proposed for deriving lower bounds for specific distance matrices. 
Among them, the branch-and-price algorithm developed by \citet{Irnich2010} and the depth-first search (DFS*) algorithm proposed by \citet{Uthus2009} have resulted in some of the best instance-specific lower bounds known to date.

\section{Approximation Results}
\label{sec:approx}

The previous section focusses on the minimal travel distance in a TTP instance.
At the same time, observe that no solution for the TTP has a distance longer than $\sum_{i\in T}\sum_{j\in T\setminus\{i\}}(d_{i,j}+d_{j,i})=2\Delta$ (i.e., each team visits every other team during a road trip of length 1). 
Combining this with \Cref{thm:maxDist}, we obtain the following result.

\begin{corollary}[{\citet{Miyashiro2012}}]
    \label{cor:approx}
	Any algorithm for the TTP(k) is a $k$-approximation algorithm.
\end{corollary}

Theoretical work on the TTP has mainly focused on improving this bound.
Since a clear overview of the best known bounds for the different variants of the TTP is currently missing in the literature, we provide it here.
For an overview of the approximation ratios, see \Cref{table:TTPApproxRegular,table:TTPApproxSpecial}. 

\subsection{TTP(2)}

From an approximation point of view, TTP(2) is by far the most studied variant.
\citet{Thielen2012} provide an algorithm with approximation ratio $1.5+\mathcal{O}(1/n)$, which is the first algorithm to achieve a ratio better than the trivial ratio of 2 from \Cref{cor:approx}.
Assuming $n/2$ is even, they provide an improved algorithm with ratio $1+16/n$.
This ratio has later been improved to about $1+4/n$ by \citet{Xiao2016}; when $n\leqslant 32$, a slightly better ratio is provided by \citet{Chatterjee2021a}.
For the case when $n/2$ is odd, \citet{Imahori2021} introduce the first $1+\mathcal{O}(1/n)$ algorithm, achieving a ratio of $1+24/n$.
Current best results are by \citet{Zhao2025a}, who attain an approximation ratio of less than $1+3/n$ and $1+5/n$ for instances where $n/2$ is even and odd, respectively.
\citet{Kanaya2025} provide a $1+9/n$ approximation algorithm, which is slightly worse than the best known one, yet it runs in $\mathcal{O}(n^3)$ rather than $\mathcal{O}(n^4)$ as required by \citet{Zhao2025a}.

\subsection{TTP(3)}

While TTP(3) poses APX-hardness, indicating the impossibility of achieving a $1+O(1/n)$ approximation ratio (unless $\mathcal{P}=\mathcal{NP}$), \citet{Fujiwara2007} demonstrate the existence of such an algorithm when restricting the problem to constant distances (i.e., CD-TTP(3)).
Their algorithm is referred to as the modified circle method.
We note that \citet{Fujiwara2007} also provide a second algorithm for CD-TTP(3), yielding optimal solutions when $n$ is a multiple of 6 minus 2 and $n\leqslant 50$.
However, this heuristic assumes the availability of a break minimum 1RR with no breaks on some pre-defined time slots (generated with IP), rendering it non-polynomial.

By randomly permuting team names and employing derandomization techniques, \citet{Miyashiro2012} transform the modified circle method into an approximation algorithm with a ratio of $2+\mathcal{O}(1/n)$ for the more general TTP(3).
Similarly, permuting teams based on a solution to the associated TSP instance, \citet{Yamaguchi2011} achieve a ratio of $1.667+\mathcal{O}(1/n)$.
\citet{Zhao2025e} utilize a 3-cycle packing to improve the ratio to $1.598+\epsilon$, for any $\epsilon>0$.

\subsection{TTP($k>3$)}

The algorithm from the previous section by \citet{Yamaguchi2011} extends to any value of $k$: for $k> 5$, the approximation ratio is $(5k-7)/(2k)+O(k/n)$, improving to $(2k-1)/k+O(k/n)$ when $k\leq 5$. 
In fact, the schedules provided by \citet{Yamaguchi2011} also satisfy the mirroring constraint, implying that these approximation ratios are also valid for the m-TTP($k$).
The first constant-factor approximation algorithm for TTP($k$), independent of both $n$ and $k$, is presented by \citet{Westphal2012}, claiming a bound of 5.875. 
However, \citet{Zhao2025d} point out that the algorithm by \citet{Westphal2012} may violate the at-most constraint and identify flaws in the analysis of the approximation ratio, suggesting the correct ratio is 6.667 rather than 5.875. By rectifying these errors and refining certain lower bounds, \citet{Zhao2025d} achieve a constant approximation ratio of 5, which further improves to 4 when $k\geqslant n/2$.
\citet{Imahori2014} achieve a ratio of 2.75 when $k = n - 1$ (i.e., the UTTP).

\subsection{Other special variants}

Regarding the linear-distance TTP, \citet{Hoshino2012} provide an expander construction method, transforming a 1RR with $n$ teams into a 2RR with $3n-2$ teams, achieving a ratio of 1.333. For multiples of 6, \citet{Zhao2022} improve the ratio to $1.2+\epsilon$, for any $\epsilon > 0$. \citet{Zhao2023a} present an EPTAS for LD-TTP($k$) with $k\geqslant 3$, theoretically providing solutions arbitrarily close to the optimum.

\citet{Hoshino2012} extend their approximation algorithm to a heuristic for TTP(3) by mapping teams on a straight line and solving the associated LD-TTP(3) instance. \citet{Rasmussen2008b} propose a similar approach, mapping teams circularly and reusing the best-known solution of the associated CD-TTP(3) instance. However, no approximation ratios are provided.

Finally, \citet{Hoshino2013} provide a $2+\mathcal{O}(1/n)$ approximation algorithm for B-TTP(3) when $n$ is a multiple of 3, which can be slightly improved further for special distance metrics like the Euclidean one.
\citet{Zhao2025b} further improve this ratio to $1.5+\epsilon$, for any $\epsilon > 0$ and any $n$.
Approximation results for a generalization of B-TTP(3) to the three-partite case can be found in \citet{Zhao2026}.

\begin{table}[t]
	\tiny
	\setlength{\tabcolsep}{5pt}	
	\centering	
	\begin{threeparttable}	
	\caption{Approximation ratios for the regular TTP}
	\begin{tabular}{lllll}
		\toprule[1pt] 
		Publications & TTP(2) & TTP(3) & TTP(4) &TTP(k)\\
		\midrule
		 \citet{Thielen2012}    & $1.5 + 6/(n-4)$\tnote{\tiny a}\\
		                        & $1.5 + 5/(n-1)$\tnote{\tiny b}\\
		                        & $1+16/n$\tnote{ \tiny b,c}\\
		 \citet{Xiao2016}       & $1+2/(n-2)+2/n$\tnote{\tiny b}\\
		 \citet{Imahori2021}    & $1+24/n$\tnote{\tiny a}\\

   	 \citet{Chatterjee2021a}& $1+ (\lceil \log_2 n/4 \rceil+4)/(2(n-2))$\tnote{\tiny b}\\
     \citet{Zhao2025a}       & $1+5/n- 10/(n(n-2))$\tnote{\tiny a}\\
		                      & $1+3/n- 10/(n(n-2))$\tnote{\tiny b}\\[10pt]
	    \citet{Miyashiro2012}  && $2+2.25/(n-1)$\\
     \citet{Yamaguchi2011}  && $1.667+O(1/n)$\\
	     	 		
		 \citet{Zhao2025e}       && $1.598 + \epsilon$ & $1.7+ \epsilon$\\
      \citet{Zhao2025c}       && & $1.625+ \epsilon$\\[10pt]
		 \citet{Yamaguchi2011}  &&&&$(2k-1)/k+O(k/n)$\tnote{\tiny d}\\
					&&&&$(5k-7)/2k+O(k/n)$\tnote{\tiny e}\\
		 \citet{Westphal2012}   &&&&5.875\tnote{\tiny f}\\
    \citet{Zhao2025d} &&&& 5.0\\
                      &&&& 4.0\tnote{\tiny g}\\
		\bottomrule[1pt]
	\end{tabular}
	\label{table:TTPApproxRegular}
	 \begin{tablenotes}
		 \setlength{\columnsep}{0.8cm}
		 \setlength{\multicolsep}{0cm}
		 \begin{multicols}{2}
		\item[a] $n/2$ is odd.
		\item[b] $n/2$ is even.
		\item[c] $n \geq 12$.
		\item[d] $k \leq 5$, also valid for m-TTP($k$).
		\item[e] $k \geq 6$, also valid for m-TTP($k$). \citet{Zhao2025e} show how to improve the ratio slightly further.
        \item[f] $k \geq 4$ and $n\geq 6$. \citet{Zhao2025d} claim the correct ratio is 6.667.
        \item[g] $k\geq n/2$.
		\end{multicols}
	\end{tablenotes}
\end{threeparttable}	
\end{table}

\begin{table}[t]
	\tiny
	\setlength{\tabcolsep}{8pt}	
	\centering	
	\begin{threeparttable}	
	\caption{Approximation ratios for special variants of the TTP}
	\begin{tabular}{lllll}
		\toprule[1pt]
		Publications & CD-TTP(3) & LD-TTP(3) & B-TTP(3) &UTTP\\
		\midrule
		\citet{Fujiwara2007}    &$1+\frac{1/3(n-1)}{4/3n^2-n}$ \tnote{\tiny a}  \\
					& $1+\frac{1/3n-1/3}{4/3n^2-5/6n-1}$ \tnote{\tiny b}\\
					&$1+\frac{5/6n-5/3}{4/3n^2-2/3n}$ \tnote{\tiny c}\\
		\citet{Hoshino2012} && 1.333\tnote{\tiny d}\\
		\citet{Zhao2022}    && $1.2+\epsilon$\tnote{\tiny e}\\
        \citet{Zhao2023a} && EPTAS\tnote{\tiny f}\\
		\citet{Hoshino2013} &&& $2+\mathcal{O}(1/n)$\tnote{\tiny a}\\	
        \citet{Zhao2025b} &&& $1.5+\epsilon$\\
	    \citet{Imahori2014}    &&&&2.75\\
		\bottomrule[1pt]
	\end{tabular}
	\label{table:TTPApproxSpecial}
	 \begin{tablenotes}
		 \setlength{\columnsep}{0.8cm}
		 \setlength{\multicolsep}{0cm}
		 \begin{multicols}{2}
		\item[a] $n\equiv 0$ (mod 3).
		\item[b] $n\equiv 1$ (mod 3).
		\item[c] $n\equiv 2$ (mod 3).
		\item[d] $n\equiv 4$ (mod 6).
  		\item[e] $n\equiv 0$ (mod 6).
        \item[f] valid for any $k\geq 3$
		\end{multicols}
	\end{tablenotes}
\end{threeparttable}	
\end{table}

\section{Optimization Algorithms}
\label{sec:optimization}

Dozens of algorithms have been proposed to solve TTP instances either optimally or heuristically.
Since there are far too many algorithms to enumerate them all, this section only focuses on a selection of them.

\subsection{Exact Methods}

The first (unsuccessful) attempt reported in the literature to solve the NL8 instance with a more sophisticated method than basic IP or CP models is by \citet{Benoist2001}, who introduce a hybrid algorithm combining Lagrangian relaxation with constraint programming. 
Shortly after, \citet{Easton2003} succeed in solving NL8 using a branch-and-price method, requiring approximately four days of computation time on 20 processors. They use CP to generate columns representing the venues at which each team plays in each round.
Later, \citet{Irnich2010} reformulate the pricing problem as a shortest path problem over an expanded network, significantly accelerating the algorithm and improving several best lower bounds at that time. 
Building upon this, \citet{Uthus2009} utilize some symmetry reductions proposed by \citet{Irnich2010} along with depth-first branch-and-bound (DFS*) to solve NL8 using only 4 processors and 100 seconds of computation time.
In a follow-up work, \citet{Uthus2012} propose an iterative-deepening algorithm (A*), achieving optimality for problem instances with up to 10 teams for the first time. However, attempts to solve larger instances encounter memory issues.
A polyhedral study for an IP formulation for the UTTP together with a new class of valid inequalities is given in \citet{Siemann2022}.

Regarding the time-relaxed TTP, \citet{Brandao2014} propose an exact approach based on combinatorial branch-and-bound combined with dynamic programming to efficiently compute the ILB. A related idea for the regular TTP is introduced by \citet{Frohner2023}, who propose a beam search heuristic that explores the most promising nodes at each level of the search tree, guided by either the ILB or a heuristic approximation thereof. This approach results in several of the currently best known solutions. 

For over a decade, apart from the constant distance instances, there has been no advancement in optimally solving instances with more than 10 teams. The current best upper bound for NL12 is 110,729, with the best lower bound being 108,629, leaving a gap of approximately 2\%. Consequently, we propose the following open problem, whose solution would constitute a significant milestone in the study of the TTP.

\begin{open}
    Determine a proven optimal solution for the NL12 instance.
\end{open}

\subsection{Heuristic Methods}

The success of the TTP may at least be partially explained by the fact that it provided one of the first optimization problems in sports timetabling that was both challenging and sufficiently simple for the application of metaheuristics, at a time when these methods were still gaining popularity.
Indeed, before the introduction of the TTP, metaheuristics for sports scheduling were almost non-existing and neighborhood structures for altering round-robin timetables were largely unexplored.

This is reflected in the early heuristics applied to this problem, which often avoided the need for neighborhood structures by combining IP or CP solvers as subroutines within a heuristic framework for achieving large moves in the search space.
For instance, \citet{Henz2004b} suggests a fix-and-optimize heuristic based on CP, where different submodels are solved that focus on the optimization for only a subset of teams or time slots, while fixing the part of the timetable that is not optimized for.
Alternatively, they fix the home-away statuses for the teams while optimizing the opponents in each round, or vice versa.
Another example is \citet{Crauwels2003} and \citet{Adriaen2003}, who use Ant Colony Optimization (ACO) based techniques, in which $n$ ants traverse a network to construct each team’s schedule, backtracking whenever infeasibilities occur.
As such, rather than relying on dedicated neighborhoods, problem-specific information is encoded in the network structure.
While the initial performance of these early ACO heuristics was somewhat disappointing, at least when viewed retrospectively compared to current state-of-the-art methods, \citet{Uthus2009} later show how the performance of this heuristic can be drastically improved by incorporating more sophisticated techniques for handling the problem constraints.

A major breakthrough comes with the introduction of two novel neighborhoods that preserve the round-robin structure, proposed by \citet{Anagnostopoulos2006} using Simulated Annealing (SA) and \citet{Ribeiro2007} using GRASP and iterated local search.
The first neighborhood, PartialRoundSwap (PRS) exchanges a subset of the games between two rounds.
The second, PartialTeamSwap (PTS) swaps the opponents of two teams in a given round. 
Both use a deterministic ejection or repair chain to restore feasibility (see \citet{Ribeiro2025} for details). 
Several follow-up works have emerged from these two pioneering papers.
For instance, \citet{VanHentenryck2006} adapt the SA algorithm to solve the mirrored TTP, \citet{Lim2006} divide the search space into a timetable and a team assignment phase, and \citet{VanHentenryck2007} propose a massively parallelized version of the SA algorithm.
\citet{Gaspero2007} incorporate the neighborhoods into a tabu search framework, while \citet{Goerigk2016} combine tabu search with integer-programming–based neighborhoods.
Neighborhoods dedicated to RTTP are proposed in \citet{Montero2015}.
Many of today’s best-known solutions stem from these approaches.

In the very beginning, little was known regarding PRS and PTS and often only restricted variants of the moves were used.
For instance, \citet{Ribeiro2007} only apply PTS moves involving 4 games, however, as shown by \citet{Gaspero2007} PTS is equivalent to PRS in this case.
\citet{Costa2012} later observe experimentally that under some circumstances none of the moves above are able to escape a given solution, and thus that the search space is disconnected (disproving a conjecture by \citet{Gaspero2007}).
These results are later formalized and proved for a wider family of timetables in \citet{Januario2015} and \citet{Januario2016c}.
A paper by \citet{Langford2010}, which to our surprise is majorly overlooked in the literature, shows how PTS in the context of 2RR tournaments can be improved through a careful selection of the SwapHomes neighborhood swapping the assignment of rounds to games $(i,j)$ and $(j,i)$.
This considerably reduces the number of opponents involved in the repair chain.
Moreover, while PRS is commonly attributed to \citet{Anagnostopoulos2006} and \citet{Ribeiro2007}, we note that it is in fact proposed more than ten years earlier by \citet{Russell1994} who apply it to a sports scheduling problem that can be considered a variant of TTP(2).

Ever since the introduction of PRS and PTS, only one new neighborhood has been proposed for time-constrained round-robin scheduling, known as TARS and its generalized variant GPTS (see \citet{Januario2016} and \citet{Ribeiro2025}).
These neighborhoods have resulted in several new best-found solutions for the TTP-PV, but have not been applied yet to the classic TTP(3).
Even with the inclusion of GPTS to the neighborhood structure, it is unknown whether the search space of round-robin tournament scheduling is connected (see \citet{Ribeiro2025}). This leads to the following open question.

\begin{open}
Does there exist a polynomial-time neighborhood that fully connects the round-robin timetable search space, while preserving the 2RR structure among consecutive moves?
\end{open}

To achieve the next major breakthrough in metaheuristic development for the TTP, we believe it would be beneficial to come up with new neighborhoods that explicitly take into account the travel cost or ensure feasibility of the no-repeater and at-most constraints (see also \citet{Caeceres2012} and \citet{Loyen2025}).

\begin{open}
\label{open:newN}
Do there exist polynomial-time neighborhoods that explicitly optimize for travel distance and/or are guaranteed to respect the no-repeater and/or the at-most constraint?
\end{open}

We believe that a good candidate for Open Problem~\ref{open:newN} is a neighborhood that focuses solely on the home-away assignment of the teams.
A first attempt in this direction has been made by \citet{Caeceres2012} who favor SwapHomes moves that result in longer yet still valid road trips.
As for the mirrored version of the problem, we believe the neighborhoods proposed in \citet{Knust2006} to modify home-away patterns in single round-robin tournaments are a good candidate.
On the negative side, \citet{Loyen2025} claim that connected neighborhoods that operate in the feasible space only, may very well not exist.
They argue that valid TTP solutions are almost maximally different from each other and thus that moving from one solution to another implies that almost all games have to be rescheduled.
We note, though, that their observations are based on a relatively small sample of timetables and that random valid solutions are almost maximally different \emph{on average}.
While this indeed suggests it is unlikely to be able to move from one feasible solution to any other in just a single (or a few) moves, it does not rule out the possibility for longer chains of moves via interconnecting solutions.
For connectivity, analyzing the number of differences to its closest neighbor could therefore be of interest.

\citet{Nakahatta2026} observe that the successful application of population-based algorithms to the TTP remains largely unexplored. 
They partially attribute this to the lack of a construction algorithm capable of generating many diverse and feasible starting solutions, ideally sampling the entire feasible solution space without bias.
The absence of such an approach is somewhat surprising, given the existence of an algorithm that results in any possible (single) round-robin tournament (see \citet{Costa2012}) and the myriad of approximation algorithms that also provide initial feasible solutions.
Moreover, \citet{Costa2012} highlights the importance of starting from solutions that differ from the canonical construction, which remains the most commonly used starting solution to date.
This brings us to the following open problem.

\begin{open}
\label{open:initial}
Does there exist an algorithm for uniformly sampling feasible solutions to TTP(3) at random?
\end{open}

\section{Conclusion and future research directions}
\label{sec:conclusion}

Over the past 25 years, the Traveling Tournament Problem has established itself as one of the central benchmark problems in sports scheduling. Its relatively simple formulation has made it an attractive testbed for a wide range of optimization techniques, while its computational complexity continues to motivate methodological advances. As this survey has shown, substantial progress has been made both theoretically and computationally. Theoretical advances are reflected in computational complexity results and approximation guarantees, while computational progress has been driven by the development of exact and heuristic algorithms that have produced increasingly high-quality solutions for many benchmark instances. Nevertheless, besides the 7 open problems mentioned in this paper, several promising research directions remain.

First, there is a need for new benchmark instances and a more standardized framework for computational evaluation. For many values of the number of teams, only a limited set of less than 10 benchmark instances is available, with most of them having a simple topological structure. Moreover, problem instances of different sizes are very similar to each other as smaller instances are contained within larger ones. All of this restricts opportunities for parameter tuning, algorithm comparison, and robustness analysis. Expanding the benchmark library and establishing more consistent standards for reporting computational performance (e.g., resources used) would facilitate fairer and more reproducible comparisons between solution approaches.

Second, the relationship between the Traveling Tournament Problem (TTP) and the Traveling Umpire Problem (TUP) deserves further investigation. We believe that a meaningful trade-off may exist between the travel distances of teams and umpires, where a modest increase in the TTP objective could yield a substantial reduction in umpire travel. To the best of our knowledge, only \citet{Bender2016} consider this integrated setting, and the trade-off between these competing objectives remains largely unexplored.

Finally, while the TTP is intentionally a simplified abstraction of real-world scheduling, this simplicity is arguably one of its greatest strengths, as it has enabled the development of a rich body of algorithmic research. Rather than replacing the TTP with increasingly complex variants, we advocate investigating how the insights and techniques developed for the TTP can be transferred to practical scheduling applications. Despite 25 years of research, only a single documented real-world application is known to us \citep{Duran2019}. Bridging the gap between theoretical advances and practical deployment therefore represents an important opportunity for future work.

We hope that this survey provides both a comprehensive overview of the existing literature and a useful starting point for researchers entering the field. We also encourage researchers to share and validate their results on the RobinX website\footnote{\url{www.robinxval.ugent.be/RobinX/}} (see \citet{VanBulck2019}), with which we continue to track the best-known solutions and bounds. Given the continued methodological developments in combinatorial optimization and the growing availability of computational resources, the Traveling Tournament Problem is likely to remain an important benchmark and source of challenging research problems for years to come.

\section*{Acknowledgements}
Fan Yang is supported by the Young Scientists Fund of the National Natural Science Foundation of China (grant number 72301177) and the Shanghai Pujiang Program (grant number 22PJC091). 
David Van Bulck is supported by the Research Foundation Flanders (FWO) [23AXE35N].

\begin{singlespace}
\renewcommand{\bibfont}{\small}
\setlength{\bibsep}{2.0pt}
\bibliography{library.bib}
\end{singlespace}


\end{document}